\documentclass{article}
\usepackage{amssymb}
\usepackage{graphicx}
\usepackage{amsmath}
\newtheorem{theorem}{Theorem}

\newtheorem{corollary}[theorem]{Corollary}

\newtheorem{definition}{Definition}[section]

\newtheorem{lemma}{Lemma}[section]

\newtheorem{problem}{Problem}
\newtheorem{proposition}{Proposition}[section]

\newcommand{\Z}{\mathbb{Z}}
\input{tcilatex}
\begin{document}

\title{Dynamical properties of profinite actions}
\author{Mikl\'{o}s Ab\'{e}rt and G\'{a}bor Elek}
\maketitle

\begin{abstract}We study profinite actions of residually finite groups in terms of weak
containment.

 We show that two strongly ergodic profinite actions of a group are weakly
equivalent if and only if they are isomorphic. This allows us to construct
continuum many pairwise weakly inequivalent free actions of a large class of
groups, including free groups and linear groups with property (T).

 We also prove that for chains of subgroups of finite index, Lubotzky's
property ($\tau$) is inherited when taking the intersection with a fixed
subgroup of finite index. That this is not true for families of subgroups in
general leads to answering the question of Lubotzky and Zuk, whether for
families of subgroups, property ($\tau$) is inherited to the lattice of
subgroups generated by the family.

 On the other hand, we show that for families of normal subgroups of finite
index, the above intersection property does hold. In fact, one can give
explicite estimates on how the spectral gap changes when passing to the
intersection.

 Our results also have an interesting graph theoretical consequence that
does not use the language of groups. Namely, we show that an expander
covering tower of finite regular graphs is either bipartite or stays bounded
away from being bipartite in the normalized edge distance.

\end{abstract}

\section{Introduction}

Let $\Gamma $ be a countable group. A measure preserving action $f$ on the
Borel probability space $(X,\mu )$ is \emph{profinite},  if there
exists a sequence of finite $\Gamma $-invariant partitions $P_{n}$ of $X$
such that $P_{n}$ consists of clopen sets, each $P_{n}$ is a refinement of $%
P_{n-1}$ and the union of $P_{n}$ generates the topology on $X$. One can
obtain all the ergodic profinite actions from the group itself as follows. A 
\emph{chain} in $\Gamma $ is a sequence $\Gamma =\Gamma _{0}\geq \Gamma
_{1}\geq \ldots $ of subgroups of finite index in $\Gamma $. Let $T=T(\Gamma
,(\Gamma _{n}))$ denote the coset tree of $\Gamma $ with respect to $(\Gamma
_{n})$ and let $\partial T$ denote the boundary of $T$. Then $\Gamma $ acts
on $\partial T$ by measure-preserving homeomorphisms; we call this action
the \emph{boundary action} of $\Gamma $ with respect to $(\Gamma _{n})$. An
especially nice case is when the chain consists of normal subgroups with
trivial intersection. Here $\partial T$ is a compact topological group,
namely the profinite completion of $\Gamma $ with respect to $(\Gamma _{n})$%
, endowed with the normalized Haar measure and $\Gamma $ maps in $\partial T$
with a dense image.

Let $f$ and $g$ be measure preserving actions of $\Gamma $ on the Borel
probability spaces $(X,\mu )$ and $(Y,\nu )$, respectively. Following \cite%
{kechbook}, we say that $f$ \emph{weakly contains} $g$ ($f\succeq g$) if for
all measurable subsets $A_{1},\ldots ,A_{n}\subseteq Y$, finite sets $%
F\subseteq \Gamma $ and $\varepsilon >0$ there exist measurable subsets $%
B_{1},\ldots ,B_{n}\subseteq X$ such that 
\begin{equation*}
\left\vert \mu (B_{i}^{\gamma }\cap B_{j})-\nu (A_{i}^{\gamma }\cap
A_{j})\right\vert <\varepsilon \text{ \ (}1\leq i,j\leq n,\gamma \in F\text{%
).}
\end{equation*}%
This means that the action $f$ can simulate $g$ with arbitrarily small error. 
A natural example for weak containment is when $g$ is
a {\it factor} of $f$, that is, when there exists a $\Gamma $-equivariant
surjective measure preserving map from $X$ to $Y$. We call $f$ and $g$ \emph{%
weakly equivalent} if $f\succeq g$ and $g\succeq f$.

We say that $f$ is \emph{strongly ergodic}, if it is ergodic and it does not
weakly contain the trivial (non-ergodic) action of $\Gamma $ on two points.

Our first theorem is a general weak containment rigidity result on strongly
ergodic actions.

\begin{theorem}
\label{contrig}Let $\Gamma $ be a countable group, let $f$ be a strongly
ergodic measure preserving action of $\Gamma $ and $g$ be a finite action of 
$\Gamma $. If $f$ weakly contains $g$ then $g$ is a factor of $f$.
\end{theorem}

When applying this to profinite actions, we get the following rigidity
result.

\begin{theorem}
\label{rigid}Let $f$ and $g$ be profinite actions of $\Gamma $ such that $f$
is strongly ergodic. If $f$ and $g$ are weakly equivalent then they are
isomorphic.
\end{theorem}

In terms of chains, isomorphism of boundary actions means that all elements
in one of the chains contains a conjugate of an element of the other chain.
This result allows one to show that a natural class of groups 
has many weakly incomparable measure preserving
actions.

\begin{theorem}
\label{manyrepr}Let $\Gamma $ be a a countable linear group with 
Kazhdan's property
(T) or a finitely generated free group.
Then $\Gamma $ has continuum many, pairwise weakly incomparable free
ergodic measure preserving actions.
\end{theorem}

The analogous question for orbit equivalence has been thoroughly
investigated in the literature. Very recently, this culminated in proving
that every countable, non-amenable group has continuously many, pairwise
orbit inequivalent free ergodic measure preserving actions (see \cite%
{epstein}). Orbit equivalence rigidity has also been investigated
specifically in the profinite case, mainly for Kazhdan groups, see the work
of Ioana \cite{ioana} and Ozawa-Popa \cite{ozapopa}.

Let $\Gamma $ be a group generated by a finite symmetric set $S$. We say
that a family of subgroups of finite index $\left\{ H_{n}\mid n\geq
1\right\} $ has property ($\tau $), if the family of Schreier graphs $%
\mathrm{Sch}(\Gamma /H_{n},S)$ forms an expander family. It is easy to see
that this property is independent of $S$.
For chains, property ($\tau $) is equivalent to saying that the boundary
action has spectral gap. While spectral gap implies strong ergodicity for
arbitrary measure preserving actions, an easy example of Schmidt \cite%
{schmidt} shows that this can not be reversed in general. However, we can
show that for boundary actions with respect to normal chains, the two
properties are in fact equivalent.

\begin{theorem}
\label{ekvtau}Let $\Gamma $ be a finitely generated group. Let $(\Gamma
_{n}) $ be a normal chain in $\Gamma $ and let $f$ denote the boundary
action of $\Gamma $ with respect to $(\Gamma _{n})$. Then $f$ is strongly
ergodic if and only if it has spectral gap.
\end{theorem}

What we actually show in this direction is that spectral gap and strong
ergodicity are equivalent for compact topological groups acted on by their
dense subgroups. Since (opposed to spectral gap, see \cite{hjorthkechris}) 
strong 
ergodicity
is an orbit equivalence invariant, we get that for these actions, having
spectral gap is an orbit equivalence invariant as well.

Our next theorem shows that Theorem \ref{ekvtau} does not hold for arbitrary
chains. Let $F_{k}$ denote the free group of rank $k$.

\begin{theorem}
\label{pelda}For every $k\geq 3$ there exists a chain $(\Gamma _{n})$ in $%
F_{k}$ such that the boundary action of $\Gamma $ with respect to $(\Gamma
_{n})$ is free and strongly ergodic but $(\Gamma _{n})$ does not have
property ($\tau $).
\end{theorem}

The proof is probabilistic; it amalgamates the random lifting method of
Friedman \cite{friedman} with results on random actions on rooted trees
investigated in \cite{abevir}.

Our next result shows that property ($\tau $) of a chain is inherited when
taking the intersection with a finite index subgroup.

\begin{theorem}
\label{tauchain}Let $\Gamma $ be a finitely generated group and let $(\Gamma
_{n})$ be a chain in $\Gamma $ with property ($\tau $). Let $H$ be a
subgroup of finite index in $\Gamma $. Then the chain $(H\cap \Gamma _{n})$
has property ($\tau $) in $H$.
\end{theorem}

This has been known for normal chains by the work of Shalom \cite{shalom}.
In the case of normal subgroups, or more generally, for compact metrizable
topological groups acted on by their dense subgroups, we obtain a stronger
result, that gives an explicite lower estimate on how the spectral gap
changes when passing to a subgroup of finite index.

\begin{theorem}
\label{nagytetel}Let $G$ be a compact metrizable topological group endowed
with its normalized Haar measure $\mu $ and let $\Gamma $ be a dense
subgroup in $G$, generated by the finite symmetric set $S$. Let $H$ be a
subgroup of $\Gamma $ of index $k$, let $C$ be a coset representative system
for $H$ in $\Gamma $ and let $T=N(S,C)$ be the Nielsen-Schreier generating
set of $H$ with respect to $S$ and $C$. Let $O$ be an ergodic component of $%
G $ under the action of $H$. Then we have 
\begin{equation*}
\mathrm{h}(O,T)>\frac{1}{8k^{3-\log _{2}3}}\min \left\{ \frac{\mathrm{h}(G,S)%
}{k^{2}},1\right\}
\end{equation*}%
where $\mathrm{h}(X,S)$ denotes the Cheeger constant of the space $X$ with
respect to the set of maps $S$.
\end{theorem}

This in turn implies the following for arbitrary \emph{families} of normal
subgroups.

\begin{theorem}
\label{normalfamily}Let $\Gamma $ be a finitely generated group and let $%
\left\{ H_{n}\mid n\geq 1\right\} $ be a family of normal subgroups of
finite index in $\Gamma $ with property ($\tau $). Let $H$ be a subgroup of
finite index in $\Gamma $. Then the family $\left\{ H\cap H_{n}\mid n\geq
1\right\} $ has property ($\tau $) in $H$.
\end{theorem}

On the other hand, as we show in Section \ref{lubzuk}, Theorem \ref%
{normalfamily} fails for a general family of subgroups of finite index.
Together with Theorem \ref{tauchain} this can be used to answer a question
of Lubotzky and Zuk \cite[Question 1.14]{lubozuk}. They asked whether if $%
\left\{ H_{n}\mid n\geq 1\right\} $ is a family of finite index subgroups in 
$\Gamma $ with property ($\tau $), then the set 
\begin{equation*}
\mathcal{L}\left( \{H_{n}\}\right) =\left\{
\bigcap_{j=1}^{k}H_{n_{j}}^{g_{j}}\mid n_{j}\in \mathbb{N}\text{, }g_{j}\in
\Gamma \right\}
\end{equation*}%
also has property ($\tau $) (note that we denote a subgroup $gHg^{-1}$ by
$H^g$). The answer is negative.

\begin{corollary}
\label{lubtau}There exists a family of finite index subgroups $\left\{
H_{n}\mid n\geq 1\right\} $ in $F_{4}$, such that $\{H_{n}\}$ has property $%
(\tau )$, but the chain $\Gamma _{n}=\cap _{k=1}^{n}H_{k}$ does not.
\end{corollary}

The counterexample family $\left\{ H_{n}\mid n\geq 1\right\} $ can be
explicitely constructed. Note, however, that because of Theorem \ref%
{normalfamily}, we do not have a negative answer for the question of
Lubotzky and Zuk if we restrict our attention to normal subgroups; so for
that case, the question is still open. \bigskip

One can exploit Theorem \ref{contrig} to obtain a purely graph theoretical
result as well. By a \emph{covering tower} of graphs, we mean a sequence $%
G_{n}$ of graphs such that for all $n\geq 1$ there is a covering map from $%
G_{n+1}$ to $G_{n}$.

\begin{theorem}
\label{coverdichotomy}Let $G_{n}$ be an expanding covering tower of $k$%
-regular graphs. Then exactly one of the following holds: \newline
1) all but finitely many of the $G_{n}$ are bipartite; \newline
2) there exists $r>0$ such that for all $n$, one needs to erase at least $%
r\left\vert G_{n}\right\vert $ edges of $G_{n}$ to make it bipartite.
\end{theorem}

Equivalently to 2), the so-called independence ratio of $G_{n}$ is bounded
away from $1/2$.

In spectral language, Theorem \ref{coverdichotomy} takes the following
equivalent form: Let $G_{n}$ be a covering tower of non-bipartite $k$%
-regular graphs. If $\lambda _{1}(G_{n})$ is bounded away from $k$ then $%
\lambda _{-}(G_{n})$ is bounded away from $-k$. Here $\lambda _{1}$ denotes
the first nontrivial eigenvalue and $\lambda _{-}$ the last eigenvalue in
order. Trivially, these results are far from being true for an arbitrary
expander family of $k$-regular graphs.

It would be interesting to see whether Theorem \ref{coverdichotomy} holds
for higher chromatic numbers as well.

\begin{problem}
Let $G_{n}$ be an expanding covering tower of $k$-regular graphs such that $%
G_{n}$ can not be legally colored by $c$ colors ($n\geq 0$). Is it true that
there exists $r>0$ such for all $n$ and all $c$-colorings of $G_{n}$, the
number of unicolored edges in $G_{n}$ is at least $r\left\vert
G_{n}\right\vert $?
\end{problem}

As one would expect, almost covers of chains in amenable groups behave quite
differently from groups having a chain with property ($\tau $). Indeed, any
two free ergodic actions of an amenable group are weakly equivalent (see
\cite{kechbook} and \cite{forweiss}),
which implies that any free boundary action of a residually finite amenable
group $\Gamma $ weakly contains any finite action of $\Gamma $. This in turn
enables us to show that every $d$-generated finite solvable group can be
simulated by a $d$-generated finite $p$-group in terms of weak containment.

\begin{theorem}
\label{frattini}Let $p$ be a prime and let $F$ be a finitely generated free
group. Then the action of $F$ on its pro-$p$ completion is weakly equivalent
to the action of $F$ on its pro-(finite solvable) completion.
\end{theorem}

Trivially, the pro-$p$ completion is a factor of the pro-solvable
completion, but the other direction is somewhat surprising. We suspect that
the same result holds for the whole profinite completion.

\bigskip

The paper is organized as follows. In Section \ref{prelim} we introduce our
notions and state some of the results used later. In Section \ref{ergsec} we
prove some general ergodic theoretical results needed later for profinite
actions. Section \ref{estimate} contains the proof of Theorem \ref{nagytetel}%
. In Section \ref{weakrigidsec} we establish the weak equivalence rigidity
results and prove Theorems \ref{contrig}, \ref{rigid} and \ref{manyrepr}. In
Section \ref{sectionnontauex} we construct the example in Theorem \ref{pelda}%
. Section \ref{lubzuk} contains the proof of Theorems \ref{tauchain},
Theorem \ref{normalfamily} and Corollary \ref{lubtau}. Section \ref{graphsec}
is about the graph theoretical consequences of our results, in particular,
we prove Theorem \ref{coverdichotomy} and its corollary on eigenvalues.
Finally, in Section \ref{amensec} we deal with amenable groups, prove
Theorem \ref{frattini} and show how to derive a recent result of Conley and
Kechris \cite{conkec} using our language.

\section{Preliminaries \label{prelim}}

This section contains the general notations and some lemmas that will be
used throughout the paper.\bigskip

\noindent \textbf{Profinite and boundary actions.} Let $\Gamma $ be a group
acting on the probability space $(X,\mu )$ by measure preserving
transformations. We say that this action is \emph{profinite}, if there
exists a sequence of finite $\Gamma $-invariant partitions $P_{n}$ of $X$
such that $P_{n}$ consists of clopen sets, each $P_{n}$ is a refinement of $%
P_{n-1}$ and the union of $P_{n}$ generates the topology on $X$.

Let $(\Gamma _{n})$ be a chain in $\Gamma $. Then the \emph{coset tree} $%
T=T(\Gamma ,(\Gamma _{n}))$ of $\Gamma $ with respect to $(\Gamma _{n})$ is
defined as follows. The vertex set of $T$ equals 
\begin{equation*}
T=\left\{ \Gamma _{n}g\mid n\geq 0,g\in \Gamma \right\}
\end{equation*}%
and the edge set is defined by inclusion, that is, 
\begin{equation*}
(\Gamma _{n}g,\Gamma _{m}h)\text{ is an edge in }T\text{ if }m=n+1\text{ and 
}\Gamma _{n}g\supseteq \Gamma _{m}h
\end{equation*}%
Then $T$ is a tree rooted at $\Gamma $ and every vertex of level $n$ has the
same number of children, equal to the index $\left\vert \Gamma _{n}:\Gamma
_{n+1}\right\vert $. The left actions of $\Gamma $ on the coset spaces $%
\Gamma /\Gamma _{n}$ respect the tree structure and so $\Gamma $ acts on $T$
by automorphisms.

The boundary $\partial T$ of $T$ is defined as the set of infinite rays
starting from the root. The boundary is naturally endowed with the product
topology and product measure coming from the tree. More precisely, for $%
t=\Gamma _{n}g\in T$ let us define $\mathrm{Sh}(t)\subseteq \partial T$, the 
\emph{shadow} of $t$ as 
\begin{equation*}
\mathrm{Sh}(t)=\left\{ x\in \partial T\mid t\in x\right\}
\end{equation*}%
the set of rays going through $t$. Set the base of topology on $\partial T$
to be the set of shadows and set the measure of a shadow to be 
\begin{equation*}
\mu (\mathrm{Sh}(t))=1/\left\vert \Gamma :\Gamma _{n}\right\vert .
\end{equation*}%
This turns $\partial T$ into a totally disconnected compact space with a
Borel probability measure $\mu $. The group $\Gamma $ acts ergodically on $%
\partial T$ by measure-preserving homeomorphisms; we call this action the 
\emph{boundary action of }$\Gamma $ with respect to $(\Gamma _{n})$. See 
\cite{grineksu} where these actions were first investigated in a measure
theoretic sense.

Another way to obtain boundary actions of a finitely generated group $\Gamma 
$ is to consider its profinite completion $G$. For every closed subgroup $H$
of $G$, the right coset space $G/H$ is a compact topological space with a
normalised Haar measure on which $\Gamma $ acts from the right. One can get
a chain leading to this action by using that $H$ is an intersection of open
subgroups in $G$. It will be convenient to use this notation in Section \ref%
{weakrigidsec}.

It is easy to see that a profinite action can be obtained as a boundary
action if and only if it is ergodic. \bigskip

\noindent \textbf{Cheeger constant, spectral gap and strong ergodicity.} Let 
$(X,\mu )$ be a probability space and let $S$ be a set of measure preserving
maps. Let us define the \emph{Cheeger constant} of $X$ with respect to $S$ as%
\begin{equation*}
\mathrm{h}(X,S)=\inf \left\{ \frac{\mu (AS\setminus A)}{\mu (A)}\mid
A\subseteq X\text{, }0<\mu (A)\leq 1/2\right\}
\end{equation*}%
where 
\begin{equation*}
AS=\left\{ as\mid a\in A,s\in S\right\} \text{.}
\end{equation*}%
Note that for a finite graph $G$, the Cheeger-constant of $G$ is defined
as 
$$Ch(G):=\inf_{A\subset G\,,|A|\leq \frac{1}{2} |G|} \frac{|L(A)|}{|A|}\,,$$
where $L(A)$ denotes the number of edges between $A$ and its complement.
Now let $\Gamma $ be a group acting on the probability space $(X,\mu )$ by
measure preserving transformations. We say that this action has \emph{%
spectral gap}, if the Koopman representation of $\Gamma $ on $L^{2}(X,\mu )$
does not contain weakly the trivial representation. Here we mean the
original weak containment notion for unitary representations \cite{kechbook}. 
We will use
the following equivalent definitions. A sequence $A_{n}$ of measurable
subsets of positive measure is called an I-sequence, if for all $\gamma \in
\Gamma $ we have 
\begin{equation*}
\lim_{n\rightarrow \infty }\frac{\mu (A_{n}\gamma \setminus A_{n})}{\mu
(A_{n})}=0
\end{equation*}%
Then by \cite{schmidt} the action of $\Gamma $ has spectral gap, if and only
if it has no I-sequences. Assume now that $\Gamma $ is generated by a finite
symmetric set $S$. Then by the above, the action of $\Gamma $ has spectral
gap if and only if $h(X,S)>0$.

Let $\Gamma $ act on a probability space $(X,\mu )$ by measure preserving
maps. A sequence of subsets $A_{n}\subseteq X$ is \emph{almost invariant},
if 
\begin{equation*}
\lim_{n\rightarrow \infty }\mu (A_{n}\diagdown A_{n}\gamma )=0\text{ for all 
}\gamma \in \Gamma
\end{equation*}%
The sequence is trivial, if $\lim_{n\rightarrow \infty }\mu (A_{n})(1-\mu
(A_{n}))=0$. We say that the action is \emph{strongly ergodic}, if every
almost invariant sequence is trivial.

In the paper, we will subsequently make use of the following lemma of
Schmidt \cite{schmidt}. Let $\mathrm{Id}_{\Gamma }$ denote the trivial
action of $\Gamma $ on one point and let $\frac{1}{2}\mathrm{Id}_{\Gamma }+%
\frac{1}{2}\mathrm{Id}_{\Gamma }$ denote its trivial action on two points,
both of measure $\frac{1}{2}$\textrm{. }

\begin{lemma}[Schmidt]
Let $\Gamma $ act on a probability space $(X,\mu )$ by measure preserving
maps. If the action is ergodic, but not strongly ergodic, then for all $%
\lambda \in (0,1)$ there exists an almost invariant sequence $A_{n}\subseteq
X$ such that $\mu (A_{n})=\lambda $ ($n\geq 0$). In particular, an ergodic
 action
is strongly ergodic if and only if it does not contain $\frac{1}{2}\mathrm{Id%
}_{\Gamma }+\frac{1}{2}\mathrm{Id}_{\Gamma }$ weakly.
\end{lemma}

\noindent \textbf{Schreier graphs, Cayley graphs and property (}$\tau $%
\textbf{).} Let $\Gamma $ be a group acting on the set $X$ by permutations
and let $S$ be a subset of $\Gamma $. Then we define the Schreier graph $%
\mathrm{Sch}(X,S)$ as follows: its vertex set is $X$ and for every $s\in S$, 
$x\in X$, there is an $s$-labeled edge going from $x$ to $xs$. When $S$ is
symmetric, that is, $S=S^{-1}$, we can think on $\mathrm{Sch}(X,S)$ as an
undirected graph. A special case is when $S$ generates $\Gamma $ and $%
X=\Gamma /H$, the set of right cosets for a subgroup $H$ of $\Gamma $; in
this case $\mathrm{Sch}(\Gamma /H,S)$ is connected. When moreover, $H$ is
normal, we define the Cayley graph $\mathrm{Cay}(\Gamma /H,S)=\mathrm{Sch}%
(\Gamma /H,S)$. Cayley graphs are vertex-transitive, that is, their
automorphism groups act transitively on the set of vertices.

Let $\Gamma $ be a finitely generated group. A set $\left\{ \Gamma
_{n}\right\} $ of subgroups of finite index in $\Gamma $ has \emph{%
Lubotzky's property }$\emph{(}\tau \emph{)}$ if for some finite, symmetric
generating set $S$ of $\Gamma $, the sequence of Schreier graphs $\mathrm{Sch%
}(\Gamma /\Gamma _{n},S)$ forms an \emph{expander family}, that is, there
exists $c>0$ such that 
\begin{equation*}
h(\mathrm{Sch}(\Gamma /\Gamma _{n},S))>c\text{ \ \ (}n\geq 0\text{)}
\end{equation*}%
where the measure on $\Gamma /\Gamma _{n}$ is defined to be uniform random.
For chains, property ($\tau $) can be expressed as follows.

\begin{lemma}
\label{bobospectral}Let $(\Gamma _{n})$ be a chain in $\Gamma $. Then $%
(\Gamma _{n})$ has property ($\tau $) if and only if the boundary action of $%
\Gamma $ with respect to $(\Gamma _{n})$ has spectral gap.
\end{lemma}

\noindent \textbf{Proof.} Let $S$ be a finite symmetric generating set for $%
\Gamma $ and let $T=T(\Gamma ,(\Gamma _{n}))$ be the coset tree. Since the
set of shadows generates the topology on $\partial T$, one gets that 
\begin{equation*}
h(\partial T,S)=\inf_{n\geq 0}h(\mathrm{Sch}(\Gamma /\Gamma _{n},S))
\end{equation*}

This proves the lemma. $\square $

\bigskip

\noindent \textbf{A covering lemma.} We will use the following lemma from 
\cite{abnik}. Since we cite it in modified form, we include a short proof.

\begin{lemma}
\label{lefedes}Let $G$ be a compact topological group with normalized Haar
measure $\mu $ and let $A,B\subseteq G$ be measurable subsets of positive
measure. Let $g$ be a $\mu $-random element of $G$. Then the expected value%
\begin{equation*}
E(\mu (Ag\cap B))=\mu (A)\mu (B)\text{. }
\end{equation*}%
In particular, for any natural number $k$ there exists a subset $X$ of size $%
k$ such that 
\begin{equation*}
\mu (AX)\geq 1-(1-\mu (A))^{k}
\end{equation*}%
For $k=\left\lceil 1/\mu (A)\right\rceil $ this gives 
\begin{equation*}
\mu (AX)>1-\frac{1}{e}\text{.}
\end{equation*}
\end{lemma}

\noindent \textbf{Proof.} Let $U=\left\{ (a,g)\in G\times G\mid a\in A\text{%
, }ag\in B\right\} $. Then $U$ is measurable in $G\times G$ and using
Fubini's theorem both ways, we get 
\begin{equation*}
\mu (A)\mu (B)=\int_{a\in A}\mu (a^{-1}B)=\mu ^{2}(U)=\int_{g\in G}\mu
(Ag\cap B)=E(\mu (Ag\cap B))\,\text{.}
\end{equation*}%
The equality $E(\mu (AX))=1-(1-\mu (A))^{k}$ follows by induction on $k$.
This implies both inequalities. $\square $

\section{Strong ergodicity and spectral gap for finite index subgroups \label%
{ergsec}}

This section analyzes what happens to the strong ergodicity and spectral gap
properties for a general measure preserving action when restricting it to a
subgroup of finite index.

\begin{lemma}
\label{strongfinite}Let $\Gamma $ act ergodically on a probability space $%
(X,\mu )$ by measure preserving maps. Let $H\leq \Gamma $ be a subgroup of
finite index and let $O$ be an ergodic component of $X$ for the action of $H$%
. Then $\mu (O)$ is a multiple of $1/\left\vert \Gamma :H\right\vert $ and
the action of $\Gamma $ on $X$ is strongly ergodic if and only if the action
of $H$ on $O$ is strongly ergodic.
\end{lemma}

\noindent \textbf{Proof.} Let $H^{\prime }=\left\{ \gamma \in \Gamma \mid
O\gamma =O\right\} $ be the setwise stabilizer of $O$. Let $C$ be a coset
representative system for $H^{\prime }$ in $\Gamma $. Then $OC$ is invariant
under $\Gamma $ and hence is equal to $X$. For $x\in X$ let $f(x)$ be the
number of sets in the form of $Oc$, $c\in C$ that contain $x$. Clearly, $f$
is a measurable function. On the other hand, if $\gamma\in \Gamma$, then
$f(x)=f(x\gamma)$. Indeed, if $x$ is covered by $Oc_1, Oc_2,\dots Oc_i$,
then $x\gamma$ is covered by the $Od_1, Od_2,\dots, Od_i$, where $d_j$ is the
coset representative of $H^{\prime}c_j\gamma$. Since the $\Gamma$-action is
ergodic the function $f$ is almost everywhere equals to a constant
$l$. Therefore, $l=\mu(O)|\Gamma:H^{\prime}|$. Since $H^{\prime }\geq H$,
$\mu(O)$ is a constant multiple of  $1/\left\vert \Gamma :H\right\vert $.

Assume that the action of $H$ on $O$ is not strongly ergodic but the action
of $\Gamma $ on $X$ is strongly ergodic. Let $T$ be a coset representative
system for $H$ in $\Gamma $. Then by Schmidt's lemma, there exists an almost 
$H$-invariant sequence of measurable subsets $A_{n}\subseteq O$ such that 
\begin{equation*}
\mu (A_{n})=\frac{1}{2\left\vert \Gamma :H\right\vert }\text{ \ }(n\geq 0).
\end{equation*}
For $n\geq 0$ let $B_{n}=A_{n}T$ be the union of $T$-translates of $A_{n}$.
Let $\gamma \in \Gamma $ and for $t\in T$ let $\overline{t}\in T$ such that $%
t\gamma \overline{t}^{-1}\in H$. Then we have $1/2\left\vert \Gamma
:H\right\vert \leq \mu (B_{n})\leq 1/2$ and 
\begin{equation*}
B_{n}\gamma \setminus B_{n}\subseteq\bigcup_{t\in
  T}\left(A_{n}t\gamma\setminus A_n\overline{t}\right)\subseteq
 \bigcup_{t\in T}\left( A_{n}t\gamma 
\overline{t}^{-1}\setminus A_{n}\right) \overline{t}
\end{equation*}%
hence 
\begin{equation*}
\mu (B_{n}\gamma \setminus B_{n})\leq \sum_{t\in T}\mu (A_{n}t\gamma 
\overline{t}^{-1}\setminus A_{n})
\end{equation*}%
The latter converges to zero as $n$ tends to infinity, so $B_{n}$ is a
nontrivial almost $\Gamma $-invariant sequence, a contradiction. We get that
strong ergodicity of $\Gamma $ on $X$ implies strong ergodicity of $H$ on $O$%
.

Assume that the action of $\Gamma $ on $X$ is not strongly ergodic. 
Following the proof of Schmidt's Lemma let $%
\phi \in L^{\infty }(X,\mu )$ be the weak $\ast $-limit of a subsequence of $%
\{\chi _{A_{n}}\}_{n=1}^{\infty }$. Here we use the Banach-Alaoglu Theorem
and the fact that $L^{\infty }(X,\mu )$ is the dual of the separable Banach
space $L^{1}(X,\mu )$.

We claim that the function $\phi $ is invariant under the action of $\Gamma $%
. It is enough to see that for any Borel-set $B\subseteq X$ and $\gamma \in
\Gamma $ 
\begin{equation*}
\int_{X}(\phi \circ \gamma )\chi _{B}d\mu =\int_{X}\phi \chi _{B}d\mu \,.
\end{equation*}%
The right hand side equals to $\lim_{k\rightarrow \infty }\mu (A_{n_{k}}\cap
B)$. The left hand side equals to 
\begin{equation*}
\int_{X}\phi (\chi _{B}\circ \gamma ^{-1})d\mu =\lim_{k\rightarrow \infty
}\mu (A_{n_{k}}\cap B\gamma )=\lim_{k\rightarrow \infty
}\mu (A_{n_{k}}\gamma^{-1} \cap B)
\end{equation*}%
and our claim follows from the almost invariance of $\{A_{n}\}_{n=1}^{\infty
}$.

By ergodicity, $\phi $ is the constant $1/2$-function. This means that 
\begin{equation*}
\lim_{k\rightarrow \infty }\mu (O\cap A_{n_{k}})=\mu (O)/2\,.
\end{equation*}%
But then $O\cap A_{n_k}$ is a nontrivial almost invariant sequence with
respect to $H$. We get that the action of $H$ on $O$ is not strongly
ergodic, a contradiction.

$\square $

\bigskip

Now we will show the corresponding theorem for spectral gap. First we need a
lemma showing that if we pass to a subgroup of finite index, small sets keep
expanding.

Let $\Gamma $ be a group generated by a symmetric set $S$. Let $H$ be a
subgroup of $\Gamma $ and let $C$ be a coset representative system for $H$
in $\Gamma $. Then for each $s\in S$ and $c\in C$ there exists a unique $%
p_{c,s}\in C$ satisfying $csp_{c,s}^{-1}\in H$. Let 
\begin{equation*}
N(S,C)=\left\{ csp_{c,s}^{-1}\mid s\in S,c\in C\right\}
\end{equation*}%
It is well-known that $N(S,C)$ generates $H$.

\begin{lemma}
\label{smallexpansion}Let $\Gamma $ act ergodically on a probability space $%
(X,\mu )$ by measure preserving maps. Let $S$ be a finite symmetric
generating set for $\Gamma $ and let $H$ be a subgroup of $\Gamma $ of index 
$k$. Let $C$ be a coset representative system for $H$ in $\Gamma $ and let $%
T=N(S,C)$. Then for all measurable subsets $A\subseteq X$ with $0<\mu
(A)\leq 1/2k$ we have 
\begin{equation*}
\frac{\mu (AT\setminus A)}{\mu (A)}\geq \frac{\mathrm{h}(X,S)}{k}\text{.}
\end{equation*}
\end{lemma}

\noindent \textbf{Proof.} Let $B=AC$. Then by straightforward set
manipulations we get 
\begin{eqnarray*}
BS\setminus B &=&\bigcup_{c\in C,\text{ }s\in S}Acs\setminus \bigcup_{c\in
C}Ac=\left(\bigcup_{d\in C}\bigcup_{\substack{ c\in C,\text{ }s\in S  \\ p_{c,s}=d 
}}Acs\right)\setminus \bigcup_{d\in C}Ad\subseteq \\
&\subseteq &\bigcup_{d\in C}\left( \bigcup_{\substack{ c\in C,\text{ }s\in S 
\\ p_{c,s}=d}}Acs\setminus Ad\right) =\bigcup_{d\in C}\left( \bigcup 
_{\substack{ c\in C,\text{ }s\in S  \\ p_{c,s}=d}}Acsd^{-1}\setminus
A\right) d\subseteq \\
&\subseteq &\bigcup_{d\in C}\left( AT\setminus A\right) d=\left( AT\setminus
A\right) C
\end{eqnarray*}%
which, using the definition of the Cheeger constant and $\mu (A)\leq \mu
(B)\leq 1/2$ yields 
\begin{equation*}
k\mu (AT\setminus A)\geq \mu (BS\setminus B)\geq \mathrm{h}(X,S)\mu (B)\geq 
\mathrm{h}(X,S)\mu (A)
\end{equation*}%
from which the lemma follows. $\square $

\bigskip

Now we will show that spectral gap passes to taking a finite index subgroup.

\begin{lemma}
\label{finindexspectral}Let $\Gamma $ be a finitely generated group acting
ergodically on a probability space $(X,\mu )$ by measure preserving maps.
Let $H\leq \Gamma $ be a subgroup of finite index and let $O$ be an ergodic
component of $X$ for the action of $H$. Then the action of $\Gamma $ on $X$
has spectral gap if and only if the action of $H$ on $O$ has spectral gap.
\end{lemma}

\noindent \textbf{Proof.} Let $S$ be a finite symmetric generating set for $%
\Gamma $, let $C$ be a coset representative system for $H$ in $\Gamma $ and
let $T=N(S,C)$.

Assume that the action of $\Gamma $ on $X$ has spectral gap. Then it is also
strongly ergodic. So, by Lemma \ref{strongfinite} the action of $H$ on $%
O$ is strongly ergodic. Assume it does not have spectral gap. Then there
exists an I-sequence $A_{n}$ in $O$. So for all $h\in H$ we have 
\begin{equation*}
\lim_{n\rightarrow \infty }\frac{\mu (A_{n}h\setminus A)}{\mu (A_{n})}=0
\end{equation*}%
But that implies $\lim_{n\rightarrow \infty }\mu (A_{n})=0$, otherwise a
suitable subsequence would be a nontrivial almost invariant sequence for $H$
in $O$. Now by Lemma \ref{smallexpansion} we have 
\begin{equation*}
\dsum\limits_{t\in T}\frac{\mu (A_{n}t\setminus A)}{\mu (A_{n})}\geq \frac{%
\mu (A_{n}T\setminus A)}{\mu (A_{n})}\geq \frac{\mathrm{h}(X,S)}{\left\vert
\Gamma :H\right\vert }>0
\end{equation*}%
for all large enough $n$, a contradiction.

Now assume that the action of $H$ on $O$ has spectral gap. Then it is also
strongly ergodic, so by Lemma \ref{strongfinite} the action of $\Gamma $ on $%
X$ is strongly ergodic. Assume it does not have spectral gap. Then there
exists an I-sequence $A_{n}\subseteq X$. By strong ergodicity, we have $%
\lim_{n\rightarrow \infty }\mu (A_{n})=0$. Since $OC=X$, there exists $c\in
C $ and a subsequence $B_{n}$ of $A_{n}$ such that 
\begin{equation*}
\mu (B_{n}\cap Oc)\geq \frac{1}{k}\mu (B_{n})
\end{equation*}%
But then $(O\cap B_{n}c^{-1})$ is an $I$-sequence for $H$, a contradiction.
So the action of $\Gamma $ on $X$ has spectral gap. $\square $

\bigskip

Now we can prove a general result that will lead to Theorem \ref{ekvtau}.

\begin{proposition}
\label{strongspectral}Let $G$ be a compact topological group endowed with
its normalized Haar measure $\mu $ and let $\Gamma $ be a dense subgroup in $%
G$. Then the right action of $\Gamma $ on $G$ is strongly ergodic if and
only if it has spectral gap.
\end{proposition}

\noindent \textbf{Proof.} Spectral gap implies strong ergodicity in general,
as clearly, any nontrivial almost invariant sequence is an I-sequence.
Assume that the action of $\Gamma $ on $G$ is strongly ergodic but has no
spectral gap. Then there exists an I-sequence $A_{n}\subseteq X$. Now by Lemma
\ref{lefedes}, for any $n\geq 1$, we have a subset $X_n$ of size 
$k_n=\left\lceil 1/2\mu (A_{n})\right\rceil $ such that

\begin{equation*}
\left( 1-\frac{1}{e}\right) ^{2}<\mu (X_{n}A_{n})\leq\frac{1}{2}
\end{equation*}%
On the other hand, for all $\gamma \in \Gamma $ we have 
\begin{equation*}
X_{n}A_{n}\gamma \setminus X_{n}A_{n}\subseteq \dbigcup\limits_{x\in
X_{n}}x(A_{n}\gamma \setminus A_{n})
\end{equation*}%
which yields 
\begin{equation*}
\mu (X_{n}A_{n}\gamma \setminus X_{n}A_{n})\leq k_{n}\mu (A_{n}\gamma
\setminus A_{n})\leq \frac{\mu (A_{n}\gamma \setminus A_{n})}{2\mu (A_{n})}
\end{equation*}%
for any $\gamma \in \Gamma $. Since $A_{n}$ is an I-sequence, the last
expression converges to zero in $n$. This means that $X_{n}A_{n}$ is a
nontrivial almost invariant sequence, a contradiction. We proved that the
action of $\Gamma $ has spectral gap. $\square $

\section{Distortion of the Cheeger constant for compact groups \label%
{estimate}}

In this section we prove Theorem \ref{nagytetel} by giving an explicite
estimate on how the Cheeger constant is distorted when passing to a finite
index subgroup. We start with a general lemma on finite graphs.

\begin{lemma}
\label{graflemma}Let $G=(V,E)$ be a finite, undirected connected graph. For
a function $f:V\rightarrow \lbrack 0,1]$ let $f^{\prime }:V\rightarrow
\lbrack 0,1]$ be defined as follows. For $v\in V$ let 
\begin{equation*}
f^{\prime }(v)=\max \left\{ 0,\max \left\{ f(w)-f(v)\mid (v,w)\in E\right\}
\right\}
\end{equation*}%
and let 
\begin{equation*}
F(f)=\dsum\limits_{v\in V}f^{\prime }(v)\text{.}
\end{equation*}%
Then 
\begin{equation*}
F(f)\geq \max \left\{ f(v)\mid v\in V\right\} -\min \left\{ f(v)\mid v\in
V\right\} \text{. }
\end{equation*}
\end{lemma}

\noindent \textbf{Proof.} It is easy to see that $F$ is not increasing if we
restrict $f$ on a subgraph. Hence, by taking a path between a minimal and a
maximal element, it is enough to prove the lemma for segments. We will
proceed by erasing one of the endpoints and using induction. If any end of
the segment is not minimal or maximal, then erasing it does not change the
minimum and the maximum. The same happens if both ends are maximal. Let $v$
be an endpoint where $f$ is minimal. By erasing $v$, we may increase the
minimum, but by at most $f^{\prime }(v)$. $\square $

\bigskip

\noindent \textbf{Proof of Theorem \ref{nagytetel}.} We can assume that $%
\Gamma $ acts with spectral gap on $G$, otherwise the theorem is trivial.
Let $A$ be a measurable subset of $O$ with $0<\mu (A)\leq \mu (O)/2$. We can
also assume $\mu (A)>1/2k$, otherwise we are done by Lemma \ref%
{smallexpansion}. Let $a$ and $b$ be parameters to be set later, satisfying $%
0<a<1/2k$ and $1-1/2k<b<1$. Using $k\geq 2$, this implies $2b-1\geq a$.

Let $A_{0}=A$ and for $l>0$ let us define $A_{l+1}$ as follows. If there
exists $g\in G$ such that 
\begin{equation*}
a\mu (A_{l})\leq \mu (gA_{l}\cap A_{l})\leq b\mu (A_{l})
\end{equation*}%
then let $A_{l+1}=gA_{l}\cap A_{l}$. We do this until there is no such $g\in
G$ or $\mu (A_{l+1})\leq 1/2k$. Let $t$ be the last index and let $B=A_{t}$.
Then trivially $\mu (B)>a/2k$.

\medskip

\noindent Case 1. If $\mu (B)\leq 1/2k$ then using 
\begin{equation*}
(X\cap Y)T\setminus (X\cap Y)\subseteq (XT\setminus X)\cup (YT\setminus Y)
\end{equation*}%
we get 
\begin{equation*}
\mu (BT\setminus B)\leq 2^{t}\mu (AT\setminus A)
\end{equation*}%
which gives 
\begin{equation*}
\frac{\mu (AT\setminus A)}{\mu (A)}\geq \frac{1}{2^{t}}\frac{\mu
(BT\setminus B)}{\mu (B)}\frac{\mu (B)}{\mu (A)}\geq \frac{a}{2^{t}k}\frac{%
\mu (BT\setminus B)}{\mu (B)}
\end{equation*}%
Using Lemma \ref{smallexpansion}, this yields 
\begin{equation} \label{elsoegyenlet}
\frac{\mu (AT\setminus A)}{\mu (A)}\geq \frac{a}{2^{t}k^{2}}\mathrm{h}(G,S)  
\end{equation}

\medskip

\noindent Case 2. If $\mu (B)>1/2k$ then for all $g\in G$, we have $\mu
(gB\cap B)<a\mu (B)$ or $\mu (gB\cap B)>b\mu (B)$. Let 
\begin{equation*}
K=\left\{ g\in G\mid \mu (gB\cap B)>b\mu (B)\right\}
\end{equation*}%
Let $f,g\in K$, then using 
\begin{equation*}
fgB\setminus B\subseteq f(gB\setminus B)\cup (fB\setminus B)
\end{equation*}%
and $2b-1\geq a$ we get 
\begin{equation*}
\mu (fgB\cap B)\geq (2b-1)\mu (B)\geq a\mu (B)
\end{equation*}%
This means that $fg\in K$, so $K$ is a subgroup of $G$.

We claim that $K$ is closed. Indeed, if we approximate the indicator
function of $B$ on $G$ with a continuous function $F:G\rightarrow \mathbb{R}$
in $L^{2}$ norm well enough, then for all $g\in G$, we have 
\begin{equation*}
\mu (gB\cap B)>b\mu (B)\text{ if and only if }\dint\limits_{x\in
G}F^{g}(x)F(x)d\mu \geq \frac{a+b}{2}\mu (B)
\end{equation*}%
But the integral above is a continuous function of $g$, so our claim holds.

Let $l$ be the index of $K$ in $G$. Let $g\in G$ be a random element
according to $\mu $. Then by Lemma \ref{lefedes} the expected measure 
\begin{eqnarray*}
\mu (B)^{2} &=&E(\mu (gB\cap B))=\dint\limits_{x\in K}\mu (xB\cap B)d\mu
+\dint\limits_{x\in G\setminus K}\mu (xB\cap B)d\mu \leq \\
&\leq &\frac{1}{l}\mu (B)+\frac{l-1}{l}a\mu (B)
\end{eqnarray*}%
which gives 
\begin{equation*}
l\leq \frac{1-a}{\mu (B)-a}<\frac{2k(1-a)}{1-2ka}
\end{equation*}%
In particular, $l$ is finite and hence $K$ is open.

Let $K_{1},K_{2},\ldots ,K_{n}$ be the right cosets of $K$ in $G$ that
intersect $O$ nontrivially (in measure). Then $O\subseteq \cup _{1}^{n}K_{i}$%
, so $\mu (O)\leq n/l$. Let $B_{i}=K_{i}\cap B$ and let $p_{i}=l\mu (B_{i})$
($1\leq i\leq n$). Let 
\begin{equation*}
m=\min \left\{ p_{i}\mid 1\leq i\leq n\right\} \text{ and }M=\max \left\{
p_{i}\mid 1\leq i\leq n\right\} \text{.}
\end{equation*}%
Then%
\begin{equation*}
m\leq \frac{1}{n}\dsum\limits_{i=1}^{n}p_{i}=\frac{l}{n}\mu (B)\leq \frac{1}{%
2}\frac{l\mu (O)}{n}\leq \frac{1}{2}\text{. }
\end{equation*}%
Let $g$ be a random element of $K$ according to its normalized Haar measure $%
l\mu $. Again using Lemma \ref{lefedes} the expected measure 
\begin{equation*}
\dsum\limits_{i=1}^{l}p_{i}^{2}=E(l\mu (gB\cap B))=\dint\limits_{x\in K}l\mu
(xB\cap B)dl\mu \geq bl\mu (B)=b\dsum\limits_{i=1}^{l}p_{i}
\end{equation*}%
This gives 
\begin{equation*}
M\dsum\limits_{i=1}^{l}p_{i}\geq \dsum\limits_{i=1}^{l}p_{i}^{2}\geq
b\dsum\limits_{i=1}^{l}p_{i}
\end{equation*}%
which implies $M\geq b$.

Now $H$ acts on the partition $P=(K_{1},K_{2},\ldots ,K_{n})$ transitively
from the right, since it acts ergodically on $O$. Let $W=\mathrm{Sch}(P,T)$
be the Schreier graph for this action. We will use Lemma \ref{graflemma} on $%
W$. Let $f:P\rightarrow \lbrack 0,1]$ be defined by $f(K_{i})=p_{i}$. Then
for all $i$ we have 
\begin{equation*}
K_{i}\cap (BT\setminus B)=\dbigcup\limits_{\substack{ x\in T  \\ %
K_{j}x=K_{i} }}B_{j}x\setminus B_{i}\supseteq B_{j}x\setminus B_{i}
\end{equation*}%
for all $j\leq n$ and $x\in T$ such that $K_{j}x=K_{i}$. This implies 
\begin{equation*}
l\mu (K_{i}\cap (BT\setminus B))\geq f^{\prime }(K_{i})
\end{equation*}%
and hence, using Lemma \ref{graflemma} we get 
\begin{equation*}
\mu (BT\setminus B)\geq \frac{1}{l}F(f)\geq \frac{M-m}{l}\geq \frac{2b-1}{2l}%
\text{. }
\end{equation*}%
Again, using $\mu (A)\leq 1/2$ and the upper estimate on $l$ we get 
\begin{equation}\label{masodikegyenlet}
\frac{\mu (AT\setminus A)}{\mu (A)}\geq \frac{1}{2^{t}}\frac{\mu
(BT\setminus B)}{\mu (A)}>\frac{2b-1}{2^{t}l}>\frac{(2b-1)(1-2ka)}{%
2^{t+1}k(1-a)} 
\end{equation}

\medskip

Now we summarize the two cases. First we estimate the $t$-part. If $\mu
(B)>1/2k$, then 
\begin{equation*}
\frac{1}{2k}<\mu (A_{t})<b^{t}\mu (A)\leq \frac{1}{2}b^{t}
\end{equation*}%
which implies $b^{t}>1/k$ and so $2^{t}<k^{\log _{2}1/b}$. Otherwise, $\mu
(A_{t-1})>1/2k$ which gives $2^{t-1}<k^{\log _{2}1/b}$. Let us set the
parameters as $a=1/4k$ and $b=3/4$. We get that if $\mu (B)>1/2k$, then by
(\ref{masodikegyenlet}) we have 
\begin{equation*}
\frac{\mu (AT\setminus A)}{\mu (A)}>\frac{1}{8k}\frac{1}{2^{t}}>\frac{1}{%
8k^{3-\log _{2}3}}
\end{equation*}%
and if $\mu (B)<1/2k$ then by (\ref{elsoegyenlet}) we have 
\begin{equation*}
\frac{\mu (AT\setminus A)}{\mu (A)}\geq \frac{a}{2^{t}k^{2}}\mathrm{h}(G,S)=%
\frac{1}{8k^{3}}\frac{1}{2^{t-1}}\mathrm{h}(G,S)>\frac{1}{8k^{5-\log _{2}3}}%
\mathrm{h}(G,S)
\end{equation*}%
which in general yields 
\begin{equation*}
\frac{\mu (AT\setminus A)}{\mu (A)}>\frac{1}{8k^{3-\log _{2}3}}\min \left\{ 
\frac{\mathrm{h}(G,S)}{k^{2}},1\right\} \text{.}
\end{equation*}%
The theorem is proved. $\square $

\bigskip

\noindent \textbf{Remark.} One can probably improve the exponent of $k$ in
Theorem \ref{nagytetel} with a more careful analysis.

\section{Weak containment rigidity of profinite actions \label{weakrigidsec}}

Throughout this section we fix the following notation. Let $\Gamma $ be a
countable group. Let $(X,\mu )$ be a standard Borel probability space and $%
(Y,\nu )$ be a probability space. We allow $Y$ to be finite. Let $f$ be a
strongly ergodic measure preserving action of $\Gamma $ on $(X,\mu )$ and
let $g$ be a measure preserving action of $\Gamma $ on $(Y,\nu )$.

Let us recall the definition of weak containment for measure preserving
actions. We say that $f$ \emph{weakly contains} $g$ ($f\succeq g$) if for
all measurable subsets $A_{1},\ldots ,A_{n}\subseteq Y$, finite sets $%
F\subseteq \Gamma $ and $\varepsilon >0$ there exist measurable subsets $%
B_{1},\ldots ,B_{n}\subseteq X$ such that 
\begin{equation*}
\left\vert \mu (B_{i}^{\gamma }\cap B_{j})-\nu (A_{i}^{\gamma }\cap
A_{j})\right\vert <\varepsilon \text{ \ (}1\leq i,j\leq n,\gamma \in F\text{%
).}
\end{equation*}%
We say that $f$ \emph{contains} $g$ ($f\geq g$, or $g$ is a \emph{factor} of 
$f$) if there exists a map $\Phi :X\rightarrow Y$ which is $\Gamma $%
-equivariant and $\Phi ^{-1}(\nu )=\mu $.

Note that the names 'weakly contains' and 'contains' can be somewhat
misleading for measure preserving actions. The reason they were named like
that comes from the realm of unitary representations. In fact, $f\succeq g$
implies that the Koopman representation of $f$ weakly contains the Koopman
representation of $g$ in the unitary sense, but the reverse implication does
not hold. For details, see the recent book of Kechris \cite{kechbook}.

The first lemma says that weak containment preserves strong ergodicity and
the Cheeger constant is monotonic with respect to it.

\begin{lemma}
\label{oroklodik}Let $\Gamma $ be a countable group and let $f$ and $g$ be
measure preserving actions on the spaces $(X,\mu )$ and $(Y,\nu )$,
respectively. If $f\succeq g$ and $f$ is strongly ergodic, then $g$ is
strongly ergodic as well. Also, for any finite subset $S$ of $\Gamma $, we
have 
\begin{equation*}
\mathrm{h}(X,S)\leq \mathrm{h}(Y,S)\text{.}
\end{equation*}
\end{lemma}

\noindent \textbf{Proof.} Assume $g$ is not strongly ergodic. Let $A_{n}$ be
an almost $\Gamma $-invariant sequence of measurable subsets of $Y$ such
that $\lim_{n}\nu (A_{n})=\alpha $ with $0<\alpha <1$ ($A_{n}$ can be a
fixed set if $g$ is not ergodic). Enumerate the elements of $\Gamma $ such
that $1$ is the first element and let $F_{n}$ be the set of the first $n$
elements of $\Gamma $.

Using the weak containment condition with $A_{n},$ $F_{n}$ and $\varepsilon
=1/n,$ we get that there exist measurable subsets $B_{n}\subseteq X$ such
that 
\begin{equation*}
\left\vert \mu (B_{n})-\nu (A_{n})\right\vert <1/n
\end{equation*}%
and for all $\gamma \in \Gamma ,$ for all large enough $n$ we have 
\begin{equation*}
\left\vert \mu (B_{n}^{\gamma }\cap B_{n})-\nu (A_{n}^{\gamma }\cap
A_{n})\right\vert <1/n
\end{equation*}%
This gives us that $\lim_{n}\nu (B_{n})=\alpha $ and hence $(B_{n})$ is a
nontrivial almost invariant sequence in $X$, which contradicts the strong
ergodicity of $f$. The statement on Cheeger constants follows similarly. $%
\square $

\bigskip

Now we can prove Theorem \ref{contrig}. \bigskip

\noindent \textbf{Proof of Theorem \ref{contrig}.} Let $f$ and $g$ be
measure preserving actions on the spaces $(X,\mu )$ and $(Y,\nu )$,
respectively. Assume that $Y$ is finite, $\nu (y)\neq 0$ ($y\in Y$), $f$ is
strongly ergodic and weakly contains $g$. Then by Lemma \ref{oroklodik}, $g$
is strongly ergodic, in particular, ergodic and hence transitive. Let $%
k=\left\vert Y\right\vert $, let $y\in Y$ and let $H$ be the stabilizer of $%
y $ in $\Gamma $. Let $F_{n}$ be the first $n$ elements of $H$. Using the
weak containment condition for $\{y\},F_{n}$ and $\varepsilon =1/n$, we get
that there exists a measurable $B_{n}\subseteq X$ such that 
\begin{equation*}
\left\vert \mu (B_{n})-\nu (\{y\})\right\vert =\left\vert \mu
(B_{n})-1/k\right\vert <1/n
\end{equation*}%
and for all $\gamma \in H,$ for all large enough $n$ we have 
\begin{equation*}
\left\vert \mu (B_{n}^{\gamma }\cap B_{n})-1/k\right\vert <1/n
\end{equation*}%
That is, $B_{n}$ is a nontrivial almost 
 invariant sequence  for $H$ such that \\
$\lim_{n}\mu (B_{n})=1/k$. Now let $O_{1},\ldots ,O_{m}$ be the ergodic
components of $X$ under the action of $H$. Then for all $l\leq m$, $%
B_{n}\cap O_{l}$ is an almost $H$-invariant sequence in $O_{l}$, hence by
Lemma \ref{strongfinite} it has to be trivial. Since $\mu (O_{l})$ is a
multiple of $1/k$, we get that there exists a unique component $O$ of
measure exactly $1/k$ such that $\lim_{n}(O\setminus B_{n})=0$.

Now we define the map $\Phi :X\rightarrow Y$ as follows. For $x\in X$ there
exists $\gamma \in \Gamma $ such that $x\gamma ^{-1}\in O$. Let 
\begin{equation*}
\Phi (x)=y\gamma \text{.}
\end{equation*}%
It is easy to check that $\Phi $ is well-defined, measure preserving and $%
\Gamma $-equivariant. Hence, $g$ is a factor of $f$. $\square $

\bigskip

When $f$ is the boundary action of $\Gamma $ with respect to a chain $%
(\Gamma _{n})$, then we can say more.

\begin{lemma}
\label{weakrigid}Let $f$ be the boundary action of $\Gamma $ with respect to
a chain $(\Gamma _{n})$ and let $g$ be a finite action of $\Gamma $. If $f$
is strongly ergodic and weakly contains $g$, then there exists $n$ such that 
$g$ is a factor of the action of $\Gamma $ on $\Gamma /\Gamma _{n}$. In
particular, there exists $y\in Y$ and $n\in \mathbb{N}$ such that the
stabilizer of $y$ in $\Gamma $ contains $\Gamma _{n}$.
\end{lemma}

\noindent \textbf{Proof.} By Theorem \ref{contrig} $g$ is a factor of $f$.
Let $o\in O$ and let $U_{n}$ be the $H$-orbit on the $n$-th level of $%
T(\Gamma ,(\Gamma _{n}))$ that $o$ passes through. The sets $(U_n)$ 
define a level-transitive boundary action of $H$, hence this action is ergodic
and equals to the $H$-action on $O$. Note that the measure of $%
U_{n} $ is a multiple of $1/k$ and $\mu(U_n)$ 
converges to $1/k$, so there exists $n$
such that $U_{n}$ has measure exactly $1/k$. But then the $\Gamma $%
-translates of $U_{n}$ form a $\Gamma $-invariant partition, so $g$ (which
is isomorphic to the action of $\Gamma $ on $\Gamma /H$) is a factor of the
action of $\Gamma $ on $\Gamma /\Gamma _{n}$ and for any $u\in U_{n}$, the
stabilizer of $u$ in $\Gamma $ is contained in $H$. $\square $

\bigskip

Let $(A_{n})$ and $(B_{n})$ be chains in $\Gamma $. We say that $(A_{n})$ 
\emph{dominates} $(B_{n})$ if for all $n$ there exists $k$ and $x\in \Gamma $
such that $A_{n}^{x}\supseteq B_{k}$.

\begin{lemma}
\label{izom}If $(A_{n})$ dominates $(B_{n})$, then the boundary action of $%
\Gamma $ with respect to $(A_{n})$ is a factor of the boundary action of $%
\Gamma $ with respect to $(B_{n})$. If $(B_{n})$ also dominates $(A_{n})$,
then the boundary actions are isomorphic, that is, there exists a measure
preserving $\Gamma $-equivariant homeomorphism between $\partial T(\Gamma
,(A_{n}))$ and $\partial T(\Gamma ,(B_{n}))$.
\end{lemma}

\noindent \textbf{Proof.} Let $G$ be the profinite completion of $\Gamma $ 
\cite{wilson}.
Then $G$ acts on $T(\Gamma ,(A_{n}))$ by automorphisms and on $\partial
T(\Gamma ,(A_{n}))$ transitively by measure preserving homeomorphisms. Let $%
\overline{A_{n}}$ be the closure of $A_{n}$ in $G$, let $a=(A_{n})\in
\partial T(\Gamma ,(A_{n}))$ and let $A=\cap \overline{_{n}A_{n}}$. Then $A$
equals the stabilizer of $a$ in $G$ and the action of $G$ on $\partial
T(\Gamma ,(A_{n}))$ is isomorphic to the coset space action on $G/A$. Let us
define $\overline{B_{n}}$, $b$ and $B$ similarly using the chain $(B_{n})$.

Let 
\begin{equation*}
O_{n}=\left\{ g\in G\mid \overline{A_{n}^{g}}\supseteq B\right\}
\end{equation*}%
Then $O_{n}$ is a descending chain of non-empty closed subsets in $G$, so it
has nontrivial intersection by compactness. Thus there exists
 $g\in \cap O_{n}$ such
that $A^{g}\supseteq B$. But then the map $F:G/B\rightarrow G/A$ defined by 
\begin{equation*}
F(Bx)=Ag^{-1}x\text{ \ (}x\in G\text{)}
\end{equation*}%
is measure preserving, $\Gamma $-equivariant and surjective.

Now if $(A_{n})$ and $(B_{n})$ both dominate each other, we get that $A$ can
be conjugated into $B$ and vice versa. Since both $A$ and $B$ are closed,
they must be conjugate in $G$ (since the same is true in any finite quotient
group). Hence $F$ defined above is a bijection and the lemma holds. $\square 
$

\bigskip

We are ready to prove the weak containment rigidity theorem.

\bigskip

\noindent \textbf{Proof of Theorem \ref{rigid}.} Let $f$ and $g$ be
profinite actions for $\Gamma $. Since $f$ is strongly ergodic and is weakly
equivalent to $g$, $g$ is strongly ergodic as well and they are both
boundary actions for some chain. Let $(F_{n})$ be such a chain for $f$ and $%
(G_{n})$ be such a chain for $g$. Let $g_{n}$ be the the action of $g$ on
the $n$-th level of $T(\Gamma ,(G_{n}))$. Then $g_{n}$ is a factor of $g$,
so it is also weakly contained in $f\,$, so by Lemma \ref{weakrigid} it is a
factor of $f$. We get that every $G_{n}$ contains a suitable conjugate of
some $F_{k}$. Similarly, every $F_{n}$ contains a suitable conjugate of some 
$G_{k}$. Using Lemma \ref{izom}, we obtain that the two profinite actions are
isomorphic. $\square $

\bigskip

Now we will construct many non weakly comparable free boundary actions of a
wide class of groups. The following lemma will be useful.

\begin{lemma}
\label{profin}Let $\Gamma $ be a residually finite group and let $G$ be its
profinite completion. Let $A,$ $B$ be closed normal subgroups in $G$. Then
the action of $\Gamma $ on $G/A$ is a factor of the action of $\Gamma $ on $%
G/B$ if and only if $A\supseteq B$.
\end{lemma}

\noindent \textbf{Proof. }If $A\supseteq B$ then from the above proof $G/A$
is a factor of $G/B$.

Assume $G/A$ is a factor of $G/B$. Then there exists a chain $(A_{n})$ in $%
\Gamma $ such that $A=\cap _{n}\overline{A_{n}}$. Now $G/\overline{A_{n}}%
=\Gamma /A_{n}$ is a factor of $G/B$, and since $\Gamma /A_{n}$ is finite
and $\Gamma $ is dense in $G$, $\Gamma $-equivariance translates to being a
homomorphism. We get that $\overline{A_{n}}\supseteq B$ which yields $%
A\supseteq B$. $\square $

\bigskip

Now we will start proving Theorem \ref{manyrepr}. We will need a general
lemma on product actions that weakly contain a finite action and then a
general theorem that produces many weakly incomparable free actions of a
wide class of groups.

Let $\Gamma $ be a countable group and let $f$ and $g$ be measure preserving
actions of $\Gamma $. Then by the product action $f\times g$ we mean the
following: the underlying probability space is the product of the underlying
spaces of $f$ and $g$ and the action of $\Gamma $ on this space is the
diagonal action. The following lemma is in the genre of the classical result
that the product of a weakly mixing and an ergodic $Z$-action is ergodic.

\begin{proposition}
\label{product}Let $f,g$ be measure preserving actions of the
countable group $\Gamma $. 
Assume that $f$ is strongly ergodic and profinite, $g$ is
mixing. Then $f\times g$ is ergodic.
Suppose that $f,g$ are as above and $f\times g$ is strongly ergodic containing
the finite action $h$. Then $f$ contains $h$ as well.
\end{proposition}

\noindent \textbf{Proof.} Let $f$ be a boundary action on $(X,\mu)$ associated
to the chain $(H_n)$ and let $g$ be a mixing $\Gamma$-action on $(Y,\nu)$.
Suppose that there exists an invariant subset $A$ in $X\times Y$ such that
$$0<\lambda=(\mu\times\nu) (A)<1\,.$$
Recall that the Borel sets of $X\times Y$ can be approximated in measure
by finite union of product sets and that the Borel structure of $X$ is
generated by the shadows of the $H_n$-cosets.
Hence there exists a sequence $(B_n)$ of $H_n$-cylindrical sets
such that
\begin{equation} \label{kozelites}
\lim_{n\to\infty} (\mu\times\nu)(B_n\triangle A)=0\,.
\end{equation}
Note that a $H_n$-cylindrical set is in the form of
$$B_n=\bigcup_{x\in\Gamma/ H_n} Sh(x)\times T^x_n\,,$$
where $Sh(x)$ is the shadow of the coset $x$ and $T^x_n\subset Y$ is a
Borel-set. Let $J_n\subset H_n$ be the normal core of $H_n$ that is the
intersection of all the conjugates of $H_n$. Clearly, $J_n$ stabilizes
all the cosets in $\Gamma/H_n$.
\begin{lemma}
$$\lim_{n\to\infty}
\frac{|x\in\Gamma/H_n\,\mid\frac{\lambda}{10}<\nu(T^x_n)<1-\frac{\lambda}{10}|}
{|\Gamma:H_n|}\,.$$
\end{lemma}
\noindent \textbf{Proof.} Let $\{k_i\}^\infty_{i=1}$ be a subset of $J_n$. 
By the mixing property,
$$\lim_{i\to\infty}(B_nk_i\triangle B_n)> (1-\frac{\lambda}{10})
\frac{\lambda}{10}
\frac{|x\in\Gamma/H_n\,\mid\frac{\lambda}{10}<\nu(T^x_n)<1-\frac{\lambda}{10}|}
{|\Gamma:H_n|}\,.$$
On the other hand, by (\ref{kozelites}) and the invariance of $A$
$$\lim_{n\to\infty} \sup_{\gamma\in\Gamma} (\mu\times\nu)(B_n\gamma\triangle
B_n)=0\,.$$
Thus the lemma follows. $\square$. \vskip 0.2in \noindent
Let $R_n\subset X$ be the union of the shadows of all the $x\in\Gamma/H_n$ for
which
$\nu(T^x_n)\geq 1-\frac{\lambda}{10}$. Let $Q_n\subset X$ be the union 
of the shadows for which $\nu(T^x_n)\leq \frac{\lambda}{10}$. Clearly,
$\mu(R_n)$ does not tend to zero, since the measure of $A$ is $\lambda$.
Observe that the sets $\{R_n\}^\infty_{n=1}$ form a non-trivial almost
invariant system. Indeed, it is easy to see that for any $\gamma\in \Gamma$
$$\frac{1}{2} \mu(R_n\gamma\cap Q_n)<(\mu\times \nu)(B_n\gamma\triangle
B_n)\,.$$
Since $\mu(X\backslash (R_n\cup Q_n))$ tends to $0$ as $n$ tends to $\infty$
the almost invariance of $\{R_n\}^\infty_{n=1}$ follows.
This is in contradiction with the strong ergodicity of $f$.
Therefore, $f\times g$ is ergodic.

\smallskip
Now let us turn to the second part of our proposition.
Let $h$ be a finite action on the set $Z$ and $z\in Z$.
Let $Stab_\Gamma(z)=H$. Since $h$ is a factor of $f\times g$ a $H$-ergodic
component of $f\times g$ has measure $\frac{1}{k}$.
Assume that $h$ is not a factor of $f$. By the proof of Theorem
\ref{contrig} it follows that all the $H$-ergodic component in $X$ has measure
larger than $\frac{1}{k}$. Let $O\subset X$ be a such a $H$-ergodic component.
The action on $O$ is profinite and strongly ergodic, the $H$-action on $Y$ is
mixing thus $O\times Y$ is $H$-ergodic. Therefore the $H$-ergodic components
of $X\times Y$ have measure greater than $\frac{1}{k}$, leading to a 
contradiction. $\square$

\bigskip

\begin{proposition}
\label{dense}Let $\{G_{n}\}$ be an infinite family of non-isomorphic,\\
non-Abelian finite simple groups and let $\Gamma $ be a dense subgroup of $%
G=\dprod\limits_{n=1}^{\infty }G_{n}$ such that the right action of $\Gamma $
on $G$ has spectral gap. Then $\Gamma $ has continuously many boundary
actions that are pairwise weakly incomparable. If such a $\Gamma$ has Property
$T$ then $\Gamma $ has
continuously many free ergodic measure preserving actions that are pairwise
weakly incomparable.
\end{proposition}

\noindent \textbf{Proof.} For a subset $\alpha $ of the natural numbers let 
\begin{equation*}
G_{\alpha }=\dprod\limits_{n\in \alpha }G_{n}.
\end{equation*}%
Then $G_{\alpha }$ is a continuous image of the profinite completion of $%
\Gamma $ with kernel $K_{\alpha }$; let $f_{\alpha }$ denote the profinite
action of $\Gamma $ on $G_{\alpha }$. Then $f_{\alpha }$ is a factor of the
action of $\Gamma $ on $G$. In particular, $f_{\alpha }$ is ergodic and has
spectral gap.

Let $I$ be a collection of continuously many infinite subsets of the natural
numbers such that no two contains one another. Then for $\alpha ,\beta \in I$
with $\alpha \neq \beta $ we get that $K_{\alpha }$ is not a subset of $%
K_{\beta }$. So by Lemma \ref{profin}, $f_{\alpha }$ is not a factor of $%
f_{\beta }$ and hence by Theorem \ref{rigid} $f_{\alpha }$ does not weakly
contain $f_{\beta }$.

The actions $f_{\alpha }$ are not free in general. Let as assume that 
$\Gamma$ has property $T$.  Let $b$ denote Bernoulli
action of $\Gamma $ on the product space $\{0,1\}^{\Gamma }$. We claim that
the set 
\begin{equation*}
\left\{ f_{\alpha }\times b\mid \alpha \in I\right\}
\end{equation*}%
consists of pairwise weakly incomparable free actions.
 Freeness is trivial,
since the action $b$ is free. Let $\alpha ,\beta \in I$ be distinct and
assume that $f_{\alpha }\times b\succeq f_{\beta }\times b$. Let $g_{n}$
denote the action of $f_{\beta }$ on the $n$-th level of the corresponding
coset tree. Then $g_{n}$ is a factor of $f_{\beta }\times b$, so $f_{\alpha
}\times b\succeq g_{n}$. Then, by Kazhdan's property $T$ $f_{\alpha
}\times b\succeq g_{n}$ is strongly ergodic. Hence
 we can apply Proposition \ref{product} and get that $f_{\alpha
}\succeq g_{n}$ for all $n$. But then $f_{\alpha }\succeq f_{\beta }$, a
contradiction. So the claim holds. $\square $

\bigskip

\noindent \textbf{Proof of Theorem \ref{manyrepr}.} First let $\Gamma $ be a
linear group with property (T). Then by Strong Approximation (see \cite[page
401]{lubsegal}) $\Gamma $ has infinitely many pairwise non-isomorphic
non-Abelian finite simple quotient groups. This gives a homomorphism of $%
\Gamma $ to the product of these groups, and since any subdirect product of
non-isomorphic non-Abelian finite simple groups equals their direct product,
the image of $\Gamma $ in the product is dense. By property (T), any ergodic
measure preserving action of $\Gamma $ has spectral gap, so the assumptions
of Proposition \ref{dense} hold. 

\noindent
Now let $\Gamma=SL(2,\Z)$. It is well-known that $SL(2,\Z)$ has property $\tau$
with respect to its congruence subgroup chain. The boundary action associated
to this chain is just the natural action of $SL(2,\Z)$ on the following
product of finite simple groups,
$$G=\prod SL(2,q)\,,$$
where $q$ runs through all the prime-powers except $q=2,3$.
Therefore $SL(2,\Z)$ is a dense subgroup of $G$ and the action has spectral
gap. Also, the action is free since if $g\in SL(2,\Z)$ is in the kernel
of all the maps $\pi_p:SL(2,\Z)\to SL(2,p)$, then $q$ must be the unit
element. By the previous proposition, $SL(2,\Z)$ has continuously many
pairwise non-weakly equivalent free ergodic actions.

\noindent
Now let $H\subseteq SL(2,\Z)$ be a finite index subgroup.
Observe that there exist only finitely many $q$'s for which $H$ does not
surject onto $SL(2,q)$. Indeed, the normal core $N_A$ of $A$ has finite index
and $N_A$ either surjects onto $SL(2,q)$ or in the kernel of the quotient
map $\pi_q:SL(2,\Z)\to SL(2,q)$. In the latter case, $|SL(2,\Z):N_A|\geq
SL(2,q)$. Therefore $H$ acts densely on an infinite product of simple groups.
Since the action of $SL(2,\Z)$ on this product has a spectral gap, by Lemma
\ref{finindexspectral} the action of $H$ has spectral gap as well. Therefore,
all the finite index subgroups of $SL(2,\Z)$ have continuously many
pairwise non-weakly equivalent free ergodic actions. Now we finish the proof
of the theorem by noting the well-known fact that $SL(2,\Z)$ contains all the
finitely generated free groups as subgroups of finite index. $\square$

\section{A free strongly ergodic boundary action that is not ($\protect\tau $%
) \label{sectionnontauex}}

In this section we first introduce covers and random covering towers, then
prove Theorem \ref{pelda}. Let us outline the strategy. We will construct
two infinite covering towers of graphs $G_{n}$ and $K_{n}$. The graphs $G_{n}
$ and $K_{n}$ will have the same vertex set ($n\geq 1$) and they will stay
close in the edge metric. The tower $K_{n}$ will consist of disconnected
graphs, but with a large connected component that is an expander, while the
tower $G_{n}$ will have girth tending to infinity, but it will not be an
expander family. However, using its small distance from $K_{n}$ in the edge
metric, we will conclude that big sets still expand in $G_{n}$. Hence the
corresponding boundary action for $(G_{n})$ will be strongly ergodic but not
with spectral gap.

We will find our towers by iterating two steps. In the first step, we
perform a suitable random cover of $G_{n}$ and $K_{n}$, that does not change
the spectral gap of the large component of $K_{n}$ but increases the girth
of $G_{n}$. It is important to note that we use the $\emph{same}$ cover of $%
G_{n}$ and $K_{n}$ -- this makes sense because covers can be defined using
only the vertex set. Since simple random covers do \emph{not} increase the
girth, we will use a sequence of iterated covers, that does. Friedman's
theorem will control expansion. The girth of iterated random covers has been
first analyzed in \cite{abevir}; here we use a variant that is described in 
\cite{chains}. In the second step, we kill the Cheeger constant by using a
specific gluing technique and thus obtain $G_{n+1}$ and $K_{n+1}$.\bigskip 

Let $S$ be a set of size $k$. By an $S$\emph{-labeled graph} we mean a
finite Schreier graph for the free group $F_{S}$ on the alphabet $S$. That
is, a finite directed graph where the edges are labeled by elements of $S$
in a way such that for each vertex $v$ and $s\in S$ there is a unique $s$%
-labeled directed edge leaving $v$ and another one entering $v$. We
emphasize that the label set $S$ is not symmetric, on the contrary, the
formal inverses of elements of $S$ in $F_{S}$ are not in $S$. When needed
(especially when considering the graph as a group action), we can extend the
labeling by putting a reverse edge for each $s$-labeled edge and labeling it
by $s^{-1}$. Finally, when we talk about spectra or girth (the smallest
length of a cycle), we forget the direction and the labels and consider the
undirected graph obtained this way.

Now we will define covers for $S$-labeled graphs. The only non-standard
thing here is that we define covers just for the underlying vertex set in a
way that it simultaneously extends to any $S$-labeled graph on the
set.\bigskip

\noindent \textbf{Covers, random covers and covering towers.} Let $X$ be a
finite set and let $d>1$ be an integer. Let $\mathrm{Sym}(d)=\mathrm{Sym}%
(\left\{ 1,\ldots ,d\right\} )$ be the symmetric group on $d$ points. Let 
\begin{equation*}
Y=X\times \left\{ 1,\ldots ,d\right\} .
\end{equation*}
For a map 
\begin{equation*}
f:S\times X\rightarrow \mathrm{Sym}(d)
\end{equation*}
and an $S$-labeled graph $R$ on $X$ let us define the $S$-labeled graph $%
C_{f}(R)$ on $Y$ as follows. For $x\in X$, $k\in \left\{ 1,\ldots ,d\right\} 
$ and $s\in S$ let 
\begin{equation*}
(x,k)\cdot s=(x\cdot s,k\cdot f(s,x))
\end{equation*}
Then it is easy to check that $C_{f}(R)$ is an $S$-labeled graph and the map 
\begin{equation*}
\phi :(x,k)\mapsto x
\end{equation*}
extends to a $d$-sheeted covering from $C_{f}(R)$ to $R$.

A \emph{random }$d$\emph{-cover} of $R$ is defined as $C_{f}(R)$ where $%
f:S\times X\rightarrow \mathrm{Sym}(d)$ is chosen uniformly randomly.

Let $d_{1},d_{2},\ldots $ be a sequence of natural numbers. A random $\left(
d_{1},d_{2},\ldots ,d_{n}\right) $-cover is defined recursively as follows.
For $n=1$ let it be a random $d_{1}$-cover of $R$ and for $n>1$ let it be a
random $d_{n}$-cover of a random $\left( d_{1},d_{2},\ldots ,d_{n-1}\right) $%
-cover.

\begin{theorem}
\label{egyes}Let $X$ be a finite set and $S$ be an alphabet of size $d$.
Then there exists a constant $b<d$ such that for all $\varepsilon >0$ there
exists $k>0$ and a sequence $d_{1},d_{2},\ldots d_{k}$ of natural numbers
such that the following holds. Let $R$ be an $S$-labeled graph on $X$ and
let $R^{\prime }$ be a random $\left( d_{1},d_{2},\ldots d_{k}\right) $%
-cover of $R$ with the covering map $\phi :R^{\prime }\rightarrow R$. Then
with probability at least $1-\varepsilon $ the following hold: 
\begin{equation*}
\mathrm{girth}(R^{\prime })>\mathrm{girth}(R)
\end{equation*}%
and  
\begin{equation*}
\lambda _{1}(\phi ^{-1}(G))\leq \max \left\{ \lambda _{1}(G),b\right\} 
\end{equation*}%
for all nontrivial connected components $G$ of $R$, where $\lambda _{1}$
denotes the second largest eigenvalue of the adjacency matrix. 
\end{theorem}

We will use two non-trivial results for the proof. The first one is
essentially proved in \cite{chains} using the language of random
automorphisms acting on an infinite rooted tree.

\begin{proposition}
\label{girth}Let $d_{1},d_{2},\ldots $ be a sequence of natural numbers such
that $d_{n}\geq 2$ ($n\geq 1$) and let $G$ be a finite $S$-labeled graph.
Then for all $\varepsilon >0$ there exists $k$ such that for a random $%
\left( d_{1},d_{2},\ldots ,d_{k}\right) $-cover $G^{\prime }$ of $G,$ the
probability 
\begin{equation*}
\mathcal{P}\left( \mathrm{girth}(G^{\prime })>\mathrm{girth}(G)\right)
>1-\varepsilon .
\end{equation*}
\end{proposition}

\noindent \textbf{Proof.} Let $T$ be the rooted tree such that the root has $%
\left| G\right| $ children and every vertex of level $n>0$ has $d_{n}$
children. For each $s\in S$ assign an independent random element of the
automorphism group $\mathrm{Aut}(T)$ (in Haar measure). Let $G_{n}$ be the
Schreier graph of the action of $S$ on the $n$-th level of $T$. Then $G_{n}$
is a covering tower and the following two random variables have the same
distribution for all $n$: \newline
1) a random $\left( d_{1},d_{2},\ldots ,d_{n}\right) $-cover of $G$;\newline
2) the graph $G_{n}$, conditioned on $G_{1}=G$.

Now it is proved in \cite{chains} that almost surely, the automorphisms
assigned to $S$ generate the free group $F_{S}$ and moreover, the action of $%
F_{S}$ on the boundary of $T$ is free. This is equivalent to saying that
a.s., we have 
\begin{equation*}
\mathrm{girth}(G_{n})\longrightarrow \infty
\end{equation*}%
Since $G_{1}=G$ with positive probability and the girth is non-decreasing
for any covering tower, we get that for all $K>0$ the probability 
\begin{equation*}
\mathcal{P}\left( \mathrm{girth}(G_{n})>K\right) \longrightarrow 1
\end{equation*}%
as $n\rightarrow \infty $ and so the Proposition is proved. $\square $

\bigskip

\noindent \textbf{Remark. }It is worth to note that a single random cover
does not increase the girth a.s. (as the degree of the cover tends to
infinity). Indeed a random cover of the trivial graph (a vertex with a loop)
is just a random permutation, which has a fixed point with probability
bounded away from $0$.

\bigskip 

The second result we will use for Theorem \ref{egyes} is due to Friedman 
\cite[Theorem 1.2]{friedman}. Let the finite graph $H$ cover the graph $G$.
Then trivially, all the eigenvalues of the adjacency matrix of $G$ are also
eigenvalues for $H$. These are called the old eigenvalues of the covering
map, and the rest of the eigenvalues are called the new ones. 

\begin{proposition}[Friedman]
\label{friedmanresult}Let $G$ be a fixed graph, let $\lambda _{0}$ denote
the largest eigenvalue of $G$, and let $\rho $ denote the spectral radius of
the universal cover of $G$. Let $R_{n}(G)$ denote a uniform random $n$-fold
covering of $G$. Then there exists a positive function $\alpha (n)$ where $%
\alpha (n)\rightarrow 0$ with $n\rightarrow \infty $ such that the
probability that $H\in R_{n}(G)$ has all its new eigenvalues inside the
interval 
\begin{equation*}
\left[ -\sqrt{\lambda _{0}\rho }-\alpha (n),\sqrt{\lambda _{0}\rho }+\alpha
(n)\right] 
\end{equation*}%
goes to $1$ as $n\rightarrow \infty $. 
\end{proposition}

\bigskip 

\noindent \textbf{Proof of Theorem \ref{egyes}.} For $d$-regular graphs, the
parameters in Friedman's theorem give $\lambda _{0}=d$ and $\rho =2\sqrt{d-1}
$. This gives 
\begin{equation*}
\sqrt{\lambda _{0}\rho }=\sqrt{2d\sqrt{d-1}}
\end{equation*}%
which is bounded away from $d$. Let $a=d-\sqrt{\lambda _{0}\rho }$ and let $%
b=d-a/2$. Now using Friedman's theorem, we get that there exists $d_{1}$,
such that with probability at least $1-\varepsilon /4$, a uniform random $%
d_{1}$-cover of any $S$-labeled graph $G$ on $X_{0}=X$ will have all its new
eigenvalues inside $\left[ -b,b\right] $.  Let $X_{1}=X_{0}\times \left\{
1,\ldots ,d_{1}\right\} $ be the new underlying set. Iterating this, we get
that there exists $d_{k}$, such that with probability at least $1-\varepsilon
/2^{k+1}$, a uniform random $d_{1}$-cover of any $S$-labeled graph $G$ on $%
X_{k-1}$ will have all its new eigenvalues inside $\left[ -b,b\right] $. Let 
$X_{k}=X_{k-1}\times \left\{ 1,\ldots ,d_{k}\right\} $ be the new underlying
set. 

This will give us an infinite sequence $d_{1},d_{2},\ldots $ that satisfies
the following. For an $S$-labeled graph $R$ on $X$ let $R_{k}$ be the random
$\left( d_{1},d_{2},\ldots ,d_{k}\right) $-cover of $R$. Then with
probability at least $1-\varepsilon /2$, for any such $R$, all the new
eigenvalues of any of the $R_{k}$ will be inside $\left[ -b,b\right] $. In
particular, for all connected components $G$ of $R$ we have 
\begin{equation*}
\lambda _{1}(\phi ^{-1}(G))\leq \max \left\{ \lambda _{1}(G),b\right\} \text{%
.}
\end{equation*}

Now let us use Proposition \ref{girth} with the sequence $d_{1},d_{2},\ldots $
and $\varepsilon /2$. We get that there exists $k$ such that with
probability at least $1-\varepsilon /2$, for any $S$-labeled graph $R$ on $X$%
, the $\left( d_{1},d_{2},\ldots ,d_{k}\right) $-cover $R^{\prime }$ of $R$
will have larger girth than $R$. 

Putting the two probabilities together, the theorem is proved. $\square $

\bigskip 

\noindent \textbf{Gluing step. }Let $G,P_{1},P_{2}$ be $S$-labeled graphs
with covering maps 
\begin{equation*}
\pi _{i}:P_{i}\rightarrow G\text{ (}i=1,2\text{)}
\end{equation*}%
Let $s\in S$ and $p_{1}\in P_{1},p_{2}\in P_{2}$ such that $\pi
_{1}(p_{1})=\pi _{2}(p_{2})$. Let the $S$-labeled graph $P$ be defined as
follows. First, take the disjoint union of $P_{1}$ and $P_{2}$. Let $\pi
:P\rightarrow G$ be the union of $\pi _{1}$ and $\pi _{2}$. Now erase the $s$%
-labeled edges $(p_{1},p_{1}s)$ and $(p_{2},p_{2}s)$ and glue in the $s$%
-labeled edges $(p_{1},p_{2}s)$ and $(p_{2},p_{1}s)$.

\begin{lemma}
Assume that $\mathrm{girth}(P_{i})>2$ ($i=1,2$). Then $P$ is an $S$-labeled
graph, $\pi $ is a covering map and we have 
\begin{equation*}
\mathrm{girth}(P_{i})\geq \min (\mathrm{girth}(P_{1}),\mathrm{girth}(P_{2}))
\end{equation*}%
and 
\begin{equation*}
Ch\left( P\right) \leq \frac{2}{\min (\left\vert P_{1}\right\vert
,\left\vert P_{2}\right\vert )}
\end{equation*}
\end{lemma}

\noindent \textbf{Proof.} It is easy to check that $\pi $ is a covering map.
Since $\mathrm{girth}(P_{i})>2$, both of the removed edges are in an $s$%
-labeled cycle of length at least $3$, hence, by removing the edges $P_{i}$
stays connected and the new edges make the whole $P$ connected. Thus $P$ is
an $S$-labeled graph. A cycle in $P$ either stays in one component, or by
putting back the old edges, it becomes the disjoint union of two cycles,
hence its size is at least $\mathrm{girth}(P_{1})+\mathrm{girth}(P_{2})$.
The estimate on the Cheeger constant follows trivally by considering the
partition $P_{1}\cup P_{2}$. $\square $

\vskip 0.1in
\noindent
Let $G$ and $H$ be graphs on the same vertex set $X$. Then their edit distance
is defined as
$$d_e(G,H):=\frac{|E(G)\triangle E(H)|}{|X|}\,.$$
Let $f:S\times X\to Sym(d)$ as above and $G,H$ be $S$-labeled graphs on
$X$. Then by definition,
$d_e(G,H)=d_e(C_f(G),C_f(H)$.
Recall that for a finite graph $G$ 
$$Ch(G)=\sup_{0<|A|\leq\frac{1}{2}|V(G)|\,,A\subset V(G)}
\frac{|L(A)|}{|A|}\,,$$
where $L(A)$ is the set of edges between $A$ and its complement. Clearly,
$(G_n)$ is an expander sequence if and only if
$\liminf_{n\to\infty} Ch(G_n)\,.$
It is well-known that for any $\epsilon>0$ there exists $\delta>0$ such that
if $\lambda_0(G)-\lambda_1(G)>\epsilon$ for a $d$-regular graph $G$, then
$Ch(G)>\delta$.

\noindent
Let $G_1$ be an arbitrary $d$-regular connected $F_S$-labeled 
graph ($d=2|S|$).
Let $q$ be its second largest eigenvalue.
 Now let $b$ be the constant in Theorem \ref{egyes} and let $\delta>0$
be such a number that $Ch(G)>\delta$ if $G$ is a connected, $d$-regular graph
with second largest eigenvalue not greater than $max(q,b)$.
\begin{lemma} \label{gluelemma}
There exist two covering towers of $F_S$-Schreier graphs
$$G_1\leftarrow G_2\leftarrow\dots\quad {\mbox and}\quad
K_1\leftarrow K_2\leftarrow\dots$$
such that the following properties are satisfied.
\begin{itemize}
\item $G_1=K_1$ is a connected graph.
\item $G_n$ and $K_n$ are defined on the same vertex set and $d_e(G_n,K_n)<
\frac{\delta}{50}\,.$
\item $girth(G_n)\to\infty$, $girth(K_n)\to\infty$.
\item $Ch(G_n)\to 0$.
\item If $T_1\subset K_1, T_2\subset K_2,\dots$ are the largest components,
then $Ch(T_n)>\delta$ and $\frac{|T_n|}{|X_n|}>1-\frac{\delta}{100 d}\,.$
for any $n\geq 1$.
\end{itemize}
\end{lemma}
Before proving the lemma let see how it implies Theorem \ref{pelda}.

\vskip 0.1in
\noindent
{\bf Proof} (of Theorem \ref{pelda} )\,\,
Since $Ch(G_n)\to 0$, the boundary action associated to the covering tower does
not have a spectral gap. In order to prove that the action is strongly ergodic
it is enough to see that if $A_n\subset X_n$ and
$$\frac{1}{4}|X_n|\leq |A_n| \leq \frac {1}{3} |X_n|\,$$
then there exists at least $\frac{\delta}{20} |X_n|$ edges between $A_n$ and
its complement in $G_n$.
Observe that
$$\frac{|X_n|}{10}\leq |A_n\cap T_n|\leq\frac{1}{2}|T_n|\,.$$
Hence there are at least $\frac{\delta|X_n|}{10}$ edges between $A_n\cap T_n$
and its complement in $T_n$. Since $|K_n\backslash T_n|<\frac{\delta}{100
  d}|X_n|$ and $|E(G_n)\triangle E(K_n)|<\frac{\delta}{50}|X_n|$ we obtain
that $|E(K_n)\triangle E(T_n)|<\frac{\delta}{20} |X_n|$. Therefore there are
at least $\frac{\delta}{20}|X_n|$ edges between $A_n$ and its complement in
$K_n$. $\square$

\vskip 0.1in
\noindent
{\bf Proof (of Lemma \ref{gluelemma})} 
We construct the towers inductively. Suppose
we have already constructed $G_n$ and $K_n$. The we pick some iterated
coverings $\kappa_{n+1}:M_{n+1}\to K_n$ resp. $\kappa_{n+1}:L_{n+1}\to G_n$
 of $K_n$ resp. $G_n$ (on the same vertex set using the same $Sym(d_i)$-valued
functions) such a way that
\begin{itemize}
\item $girth(M_{n+1})> girth (K_n)$
\item $girth(L_{n+1})> girth (G_n)$
\item $L_{n+1}$ is connected.
\item $\lambda_1(\kappa_{n+1}^{-1}(G))\leq \max\{\lambda_1(G),b\}$
for any connected component of $K_n$. Particularly, each
$\kappa_{n+1}^{-1}(G)$ is connected.
\end{itemize}
The existence of such construction easily follows from Proposition \ref{girth}
 and Proposition \ref{friedmanresult}.
Now we pick another coverings $\kappa'_{n+1}:M'_{n+1}\to K_n$ resp. 
$\kappa'_{n+1}:L_{n+1}\to G_n$ satisfying the same properties. Nevertheless,
we choose $M'_{n+1}$ so large that the size of the greatest component
of $M'_{n+1}$ is still larger than
$(|M_{n+1}|+ |M'_{n+1}|)(1-\frac{\delta}{100 d})\,.$
Now let $K_{n+1}$ be the disjoint union of $M_{n+1}$ and $M'_{n+1}$ and
let $G_{n+1}$ be the union of  $L_{n+1}$ and $L'_{n+1}$ glued together.
It is easy to see that if that $M'_{n+1}$ is large enough then
$d_e(G_{n+1},K_{n+1})$ is still smaller than $\frac{\delta}{50}$. $\square$

\section{Subgroups and property ($\protect\tau $) \label{lubzuk}}

Let us first outline the contents of this section using a graph theoretical
language. Let $S$ be a finite alphabet and let $G$ be an $S$-labeled graph
(see the previous section for the definition). Label the inverses of edges
by formal inverses of elements of $S$. Now fix a symmetric set of words $%
w_{1},w_{2},\ldots ,w_{n}\in F_{S}$. Then we can define a new graph $%
G^{\prime }$ on $V$, by drawing a $w_{i}$-labeled edge from $v\in V$ to $%
v\cdot w_{i}$ ($v\in V,$ $1\leq i\leq n$). This section investigates how the
expansion of $G^{\prime }$ is related to the expansion of $G$.

If $w_{1},w_{2},\ldots ,w_{n}$ generate $F$, then it is easy to see that the
graph metric on $G^{\prime }$ is bi-Lipschitz to the one on $G$ with a
bounded Lipschitz constant and hence $G^{\prime }$ is also connected and
expansion is distorted in a bounded way. When $H=\left\langle
w_{1},w_{2},\ldots ,w_{n}\right\rangle $ is a proper subgroup of finite
index in $F$, $G^{\prime }$ may or may not stay connected. We shall present
an example for a sequence of graphs $G_{n}$ where $H$ has index $2$, $%
G_{n}^{\prime }$ stays connected but expansion vanishes. Surprisingly,
however, when $G_{n}$ comes from a chain of subgroups, or a family of normal
subgroups, expansion stays bounded away from zero -- of course, in light of
the previous claim, for chains, the bound is not absolute. This directly
leads us to answering the question of Lubotzky and Zuk. \bigskip

We start with the construction of `bad' $S$-labeled graphs. As an input, we
use a ($\tau $) chain in $F_{2}$. These exist by various arguments, see 
\cite{lubozuk} and \cite{lubkonyv}.

\bigskip

\noindent \textbf{Construction of bad Schreier graphs.} Let $F_{2}$ be
generated by $x_{1}$ and $x_{2}$. Let $(H_{n})$ be a chain in $F_{2}$ with
property ($\tau $). Let $C$ denote the cyclic group of $2$ elements
 generated by $t$, let $%
\Delta =F_{2}\times C$ and let $H_{n}^{\prime }=H_{n}\times 1\leq \Delta $.
Let 
\begin{equation*}
E_{n}=\mathrm{Sch}\left( \Delta /H_{n}^{\prime },\left\{
x_{1},x_{2},t\right\} \right)
\end{equation*}%
Then $E_{n}$ is a union of two subgraphs $E_{n,1}$ and $E_{n,t}$, both
isomorphic to 
\begin{equation*}
\mathrm{Sch}(F_{2}/H_{n},\{x_{1},x_{2}\})
\end{equation*}%
plus the action of $t$, which is a perfect matching between the two
subgraphs. Now we introduce a new generator $c$ that acts on the vertex set
of $E_{n}$ as follows. Let $e_{n,1},e_{n,2}\in E_{n,1}$ and let $%
e_{n,3}=e_{n,2}\cdot t\in E_{n,t}$. Let 
\begin{equation*}
c=\left( e_{n,1},e_{n,2},e_{n,3}\right)
\end{equation*}%
be the $3$-cycle moving only these points. Let $G_{n}$ be $E_{n}$ plus the
additional $c$-edges.

Let $\Gamma $ be the free group on the generating set $\left\{
x_{1},x_{2},t,c\right\} $. Then $\Gamma $ acts transitively on $G_{n}$. Let 
\begin{equation*}
\Gamma _{n}=\mathrm{Stab}_{\Gamma }(e_{n,2})
\end{equation*}
By transitivity, we have 
\begin{equation*}
G_{n}=\mathrm{Sch}\left( \Gamma /\Gamma _{n},\left\{ x_{1},x_{2},t,c\right\}
\right)
\end{equation*}
Note that $\Gamma _{n}$ is not a chain anymore, as $c$ ruins $G_{n}$ being a
covering tower.

Let $H$ be the kernel of the projection $\phi :\Gamma \rightarrow C$ defined
by $\phi (x_{1})=\phi (x_{2})=\phi (c)=1$ and $\phi (t)=t$. Then $H$ is a
normal subgroup of $\Gamma $ of index $2$ and by the Nielsen-Schreier
theorem, it is generated by 
\begin{equation*}
T=\left\{ x_{1},x_{2},c,tx_{1}t^{-1},tx_{2}t^{-1},tct^{-1},t^{2}\right\} 
\text{.}
\end{equation*}

\begin{proposition}
\label{rossztau}The family $\Gamma _{n}$ has property ($\tau $) in $\Gamma $
but the family $H\cap \Gamma _{n}$ does not have property ($\tau $) in $H$.
\end{proposition}

\noindent \textbf{Proof.} The sequence $E_{n}$ is an expander family, hence $%
G_{n}$ is an expander family as well. Thus the family $\Gamma _{n}$ has
property ($\tau $) in $\Gamma $.

For all $n\geq 1$ the element $ct^{-1}$ fixes $e_{n,2}$, so $ct^{-1}\in
\Gamma _{n}$ but $ct^{-1}\notin H$. This shows that $\Gamma _{n}\nleq H$.
Since $H$ has index $2$ in $\Gamma $, we get $\Gamma _{n}H=\Gamma $ and so $%
\mathrm{Sch}\left( H/H\cap \Gamma _{n},T\right) $ is isomorphic to $\mathrm{%
Sch}\left( \Gamma /\Gamma _{n},T\right) $. Moreover, in this action we have 
\begin{equation*}
x_{1}=tx_{1}t^{-1}\text{, }x_{2}=tx_{2}t^{-1}\text{ and }t^{2}=1.
\end{equation*}
Let us look at the set $E_{n,1}\subseteq \Gamma /\Gamma _{n}$. Both $x_{1}$
and $x_{2}$ fix $E_{n,1}$ as a set, so there are exactly $4$ edges in $%
\mathrm{Sch}\left( \Gamma /\Gamma _{n},T\right) $ that leave $E_{n,1}$, that
is, the edges coming from $c$ and $tct^{-1}$. This implies that the Cheeger
constant 
\begin{equation*}
\mathrm{Ch}\left( \mathrm{Sch}\left( H/H\cap \Gamma _{n},T\right) \right)
\leq \frac{4}{\left| E_{n,1}\right| }=\frac{2}{\left| \Gamma :\Gamma
_{n}\right| }.
\end{equation*}
Hence, $\mathrm{Sch}\left( H/H\cap \Gamma _{n},T\right) $ is not an
expander family in $H$ and so the family $H\cap \Gamma _{n}$ does not have
property ($\tau $) in $H$. $\square $

\bigskip

We are ready to prove Theorem \ref{tauchain}. Note that we are not aware of
any proof that does not use compactness in some form; the fact that there
are no bounds on how bad expansion can be distorted makes it dubious that
such proof exists. Even for normal chains, where by Theorem \ref{nagytetel}
there is an explicit lower bound on distortion, the only other proof we know 
\cite{shalom} uses invariant means. \bigskip

\noindent \textbf{Proof of Theorem \ref{tauchain}.} Let $S$ be a finite
symmetric generating set for $\Gamma $, let $k=\left\vert \Gamma
:H\right\vert $. Let $T=T(\Gamma ,(\Gamma _{n}))$ be the coset tree and
let $t_{n}\in T$ be the vertex representing the subgroup $\Gamma_n$. 
Then $\{t_{n}\}$ forms a ray in $T$, since $%
(\Gamma _{n})$ is a chain. Let $t\in \partial T$ denote this ray as a
boundary point. Let $T_{n}$ denote the $n$-th level of $T$. Let $O_{n}$ be
the orbit of $t_{n}$ in $T_{n}$ under the action of $H$. Then the
permutation action of $H$ on $O_{n}$ is isomorphic to the coset action of $H$
on $H/H\cap \Gamma _{n}$. Also, the union of $O_{n}$ forms a subtree, that
is isomorphic to the coset tree $T(H,(H\cap \Gamma _{n}))$ and the limit of
the $O_{n}$ equals the ergodic component of $\partial T$ under the action of 
$H$ that contains $t$. Let us call this component $O$. Now $(\Gamma _{n})$
has property ($\tau $) in $\Gamma $, so by Lemma \ref{bobospectral}, the
action of $\Gamma $ on $\partial T$ has spectral gap. Now using Lemma \ref%
{finindexspectral}, we get that the action of $H$ on $O$ also has spectral
gap. But the action of $H$ on $O$ is isomorphic to the boundary action of $H$
with respect to $(H\cap \Gamma _{n})$, so again by Lemma \ref{bobospectral}, 
$(H\cap \Gamma _{n})$ has property ($\tau $) in $H$. $\square $

\bigskip

Theorem \ref{tauchain} has been proved for normal chains by Shalom \cite%
{shalom} using invariant means. Theorem \ref{nagytetel} allows us to extend
his result to arbitrary families of normal subgroups as stated in Theorem %
\ref{normalfamily}.

\bigskip

\noindent \textbf{Proof of Theorem \ref{normalfamily}.} Let $S$ be a finite
symmetric generating set for $\Gamma $, let $k=\left\vert \Gamma
:H\right\vert $, let $C$ be a coset representative system for $H$ in $\Gamma 
$ and let $T=N(S,C)$. Let $G_{n}=\Gamma /\Gamma _{n}$. Then $G_{n}$ is a
compact (in fact, finite) topological group and the image of $\Gamma $ in $%
G_{n}$ is dense (being equal to $G_{n}$). Let $O_n$ be the orbit of $H$ in $%
G_{n}$ containing the identity of $G_{n}$. Then we can invoke Theorem \ref%
{nagytetel} and get that 
\begin{equation*}
\mathrm{h}(O_n,T)>\frac{1}{8k^{3-\log _{2}3}}\min \left\{ \frac{\mathrm{h}%
(G_{n},S)}{k^{2}},1\right\}.
\end{equation*}%
In particular, since $\mathrm{h}(G_{n},S)$ is bounded below and $k$ is
fixed, the family of Cayley graphs $\mathrm{Cay}(H/H\cap \Gamma _{n},T)$ ($%
n\geq 1$) is an expander family and so the theorem holds. $\square $

\bigskip

Finally, Proposition \ref{rossztau} and Theorem \ref{tauchain} together
allow us to answer the question of Lubotzky and Zuk.\bigskip

\noindent \textbf{Proof of Corollary \ref{lubtau}.} Let $\Gamma =F_{4}$ and
let $H\leq \Gamma $ and $\Gamma _{n}\leq \Gamma $ ($n\geq 1$) be defined as
in the construction above. Then by Proposition \ref{rossztau} the family $%
\Gamma _{n}$ has property ($\tau $). For $n\geq 1$ let $H_{n}=\cap
_{i=1}^{n}\Gamma _{i}$. Assume that the chain $(H_{n})$ has property ($\tau $%
). Then using Theorem \ref{tauchain}, the chain $(H\cap H_{n})$ has property
($\tau $) in $H$. But $H\cap H_{n}\leq H\cap \Gamma _{n}$ ($n\geq 1$) which
implies that the family $H\cap \Gamma _{n}$ ($n\geq 1$) also has property ($%
\tau $) in $H$. This contradicts Proposition \ref{rossztau}.

Hence, the chain $(H_{n})$ does not have property ($\tau $) and the
corollary is proved. $\square $

\section{Almost covers of graphs and the distance from being bipartite \label%
{graphsec}}

For general unlabeled graphs, weak containment translates as follows.

\begin{definition}
Let $G$ and $H$ be finite $k$-regular graphs. A map $f:E(G)\rightarrow E(H)$
is an $\varepsilon $\emph{-covering}, if $f$ is surjective and there exists $%
X\subseteq V(G)$ with $\left\vert X\right\vert >(1-\varepsilon )\left\vert
V(G)\right\vert $ such that for all $x\in X$ there exists $y\in V(H)$ such
that $f$ is a bijection between the set of edges leaving $x$ and the set of
edges leaving $y$.
\end{definition}

That is, $f$ is a local isomorphism at most vertices of $G$. Note that for $%
\varepsilon =0$ we get back the original notion of a finite sheeted covering
map and by our definition, $y$ is a unique function of $x$, that is, $f$
induces a map $V(G)\rightarrow V(H)$. It is easy to see that if $H$ is
connected, then every vertex in $H$ has the same number of preimages.

A sequence of finite graphs $(G_{n})$ \emph{almost covers} a finite graph $H$
if for all $\varepsilon >0$ there exists $n_{0}$ such that for all $n>n_{0}$%
, $G_{n}$ has an $\varepsilon $-covering to $H$.

By a \emph{covering tower of graphs}, we mean a sequence $(G_{n},f_{n})$ of
graphs and maps such that for all $n\geq 1$, $f_{n}$ is a covering map from $%
G_{n+1}$ to $G_{n}$. Let $(G_{n},f_{n})$ be a covering tower of connected $k$%
-regular graphs. Then we define the \emph{covering tree} $T=T(G_{n})$ as
follows. Let the vertex set of $T$ be the disjoint union of the $V(G_{n})$
and for all $n>1$ and $x\in V(G_{n})$ connect $x$ to its image under the
covering map. Then $T$ is a spherically homogeneous rooted tree. Let $%
\partial (G_{n})=\partial T$ be the boundary of the tree, that is, the set
of infinite rays in $T$, endowed with the product topology and measure. The
boundary $\partial T$ is naturally endowed with a graph structure: we
connect $(x_{n}),(y_{n})\in \partial T$ if $x_{n}$ and $y_{n}$ are connected
in $G_{n}$ for every $n$. This gives us a $k$-regular graphing, that we call
the \emph{boundary graphing} of $(G_{n},f_{n})$ and denote it by $\partial
(G_{n},f_{n})$. By composing covering maps and taking a limit, we get a
continuous covering map from the boundary graphing to $G_{n}$.

We are ready to prove Theorem \ref{coverdichotomy} after a lemma that is
folklore in graph theory.

\begin{lemma}
\label{biplemma}Let $G$ be a finite undirected $k$-regular graph and let $S$
be an alphabet on $k$ letters. Then $G$ can be turned into an $S$-labeled
graph such that every edge of $G$ is used exactly once in each direction.
\end{lemma}

\noindent \textbf{Proof.} Let $A\ $be the adjacency matrix of $G$. Then let us
look at $A$ as the adjacency matrix of a bipartite graph obtained by doubling
the vertices of $G$. It is $k$%
-regular, so it is a disjoint union of $k$ perfect matchings. That is, $A$
is the sum of $k$ permutation matrices. Let us label the directed edges of $%
G $ according to these permutations. This gives the required decomposition. $%
\square $\bigskip

Using Lemma \ref{biplemma} and putting in formal inverses of elements of $S$%
, one can turn a $k$-regular graph $G$ to a Schreier graph for $F_{S}$, such
that each edge is used exactly twice by the generators and its inverses (in
each direction). Note that if the directed edge $(x,y)$ is labeled by $s$ and 
$(y,x)$ is
labeled by $t$, then for the associated $F_S$-action $xs=y$ and $xt^{-1}=y$.

\bigskip

\noindent \textbf{Proof of Theorem \ref{coverdichotomy}.} Let $(G_n)$ be an 
expanding covering tower of graphs. Consider the
associated $F_S$-action on $G_1$ and pull back the action onto all the
covering graphs. Then we obtain a $F_S$-chain with boundary action on
$\partial T$. Let us consider the homomorphism $\phi:F_S\to C=\{1,t\}$, where
$\phi(s)=t$ for all the generators. Let $H$ be the kernel of $\phi$ a subgroup
of index $2$. Observe that the $F_S$ action $f$ on $\partial T$ has a spectral
gap since $(G_n)$ is an expander system. That is $f$ is strongly ergodic.
Now consider $\partial T$ as a $H$-space.

\vskip 0.1in
\noindent
{\bf Case 1.}\,\,\, Suppose that the $H$-action on $\partial T$ is ergodic.
Then by Lemma \ref{finindexspectral} it has a spectral gap.  Let
$g$ be the $F_S$-action on the set $\{1,t\}$ induced by $\phi$. By the
ergodicity assumption, $g$ is not a factor of $f$. Hence by Theorem
\ref{contrig}
, $f$
does not contain $g$ weakly. Let $r_n|V(G_n)|$ be the minimal number of edges
one needs to erase to make $G_n$ bipartite (with partition sets $A_n$,$B_n$).
Clearly, $r_1\geq r_2\geq\dots$
Suppose that $\lim_{n\to\infty} r_n=0$. Let $C_n$ be the shadow of $A_n$ and
  $D_n$ be the shadow of $B_n$. It is easy to see that
$\mu(C_n)\to 1/2,\, \mu(D_n)\to 1/2$ and for any $\gamma\in H$
$\mu(C_n\gamma\cap C_n)\to 0$ and $\mu(D_n\gamma\cap D_n)\to 0$.
Hence $f$ weakly contains $g$ leading to a contradiction. Therefore
 $\lim_{n\to\infty} r_n>0$.
\vskip 0.1in
\noindent
{\bf Case 1.}\,\,\, There exists a $H$-ergodic component $O$ of size $1/2$.
Similarly as in Lemma \ref{weakrigid}, this implies that if $n$
is larger than some constant $n_k$ there are exactly two $H$-orbits on the
$n$-th level. That is $G_n$ is bipartite if $n>n_k$. $\square $

\bigskip

Theorem \ref{contrig} suggests the following problem.

\begin{problem}
Let $(G_{n})$ be an expanding covering tower of $k$-regular graphs and $H$ a
finite graph such that $(G_{n})$ almost covers $H$. Does it follow that
there exists $n$ such that $G_{n}$ covers $H$?
\end{problem}

By Theorem \ref{coverdichotomy} the answer is affirmative when $H$ is a
graph with two points and $k$ edges going between them. \bigskip

On spectral language, Theorem \ref{coverdichotomy} takes the following
equivalent form. For a $k$-regular undirected graph $G$ on $v$ vertices let $%
\lambda _{0}(G)\geq \lambda _{1}(G)\geq \ldots \geq \lambda
_{v-1}(G)=\lambda _{-}(G)$ denote the eigenvalues of the adjacency matrix of 
$G$. Then $\lambda _{0}(G)=k$ and $\lambda _{-}(G)\geq -k$. Assuming that $G$
is connected, $\lambda _{-}(G)=-k$ if and only if $G$ is bipartite.

\begin{corollary}
\label{spectraldichotomy}Let $(G_{n})$ be a covering tower of non-bipartite $k$%
-regular graphs. If $\lambda _{1}(G_{n})$ is bounded away from $k$ then $%
\lambda _{-}(G_{n})$ is bounded away from $-k$.
\end{corollary}

\noindent \textbf{Proof. }Let $G(V,E)$ be a finite $d$-regular connected
 graph. If $S$ is a subset of $V$ then let
$e(S)$ be the minimal number of edges to be removed from the graph
spanned by $S$ to make it bipartite. Let $k(S)$ be the number of edges
to be removed to disconnect $S$ from $V-S$.
Let $r(G):=\frac{e(G)}{|V|}$ and $c(G):=\min_{S\subset V\,,\, |S|\leq
  \frac{1}{2}|V|} \frac {k(S)}{|S|}$.
Desai and Rao introduced the following constant :
$$\psi(G)=\min_{S\subset V} \frac{e(S)+k(S)}{|S|}\,$$
and proved (Theorem 3.2)
that for the smallest eigenvalue of the adjacency matrix of $G$,
$q_n(G)$ 
$$q_n(G)\geq -d+\frac{\psi^2(G)}{4d}\,.$$
\begin{lemma} We have
$$\psi(G)\geq \min\{c(G),\frac{r(G)c(G)}{2d},\frac{r(G)}{4}\}$$
\end{lemma}
{\bf Proof:}\,\,
Let $w(S)=(e(S)+k(S))/|S|$.
First let $|S|\leq |V|/2 $ then $c(G)\leq w(S)$.
Now let $|V|/2 \leq |S| \leq (1-r(G)/2d)|V|$.
Then $k(S)\geq r(G)c(G)|V|/2d$, that is $w(S)\geq
r(G)c(G)/2d$.
Finally, let $(1-r(G)/2d)|V|\leq |S| \leq |V|$.
Then the number of edges in the span of $V-S$ is at most $r(G)|V|/4$
and the number of edges between $S$ and $V(S)$ is at most $r(G)|V|/2$.
Hence in order to make $S$ bipartite one needs to remove at least
$r(G)|V|/4$ edges. Otherwise, one can make $G$ bipartite by removing
less than $e(G)$ edges. Consequently, $w(S)\geq r(G)/4$.
This ends the proof of our lemma.\,\, $\square$

\bigskip

Trivially, all these results are far from being true for an arbitrary
expander sequence of $k$-regular graphs. \bigskip

\noindent \textbf{Remark. }A standard example for a sequence of finite $k$%
-regular graphs where the girth (the minimal size of a cycle) tends to
infinity and the independence ratio is bounded away from $1/2$ is due to
Bollobas \cite{bollobas} who showed that large random $k$-regular graphs
satisfy these properties. Now Theorem \ref{coverdichotomy} allows us to find
these sequences in abundance. Indeed, take the free product $\Gamma =\mathbb{%
Z}/2\mathbb{Z}\ast \mathbb{Z}/2\mathbb{Z\ast \cdots \ast Z}/2\mathbb{Z}$
(with $k$ factors), or alternatively, for an even $k\geq 4$, the free group $%
\Gamma =F_{k/2}$. Let $S$ be a standard generating set of $\Gamma $ and let $%
N$ be the kernel of the homomorphism $\Gamma \rightarrow \mathbb{Z}/2\mathbb{%
Z}$ that sends all elements of $S$ to the nontrivial element. Then by
Theorem \ref{coverdichotomy}, for any chain $(\Gamma _{n})$ in $\Gamma $
which is property $(\tau )$ and satisfies $\Gamma _{n}\nleq N$ for all $n$,
the sequence of Schreier graphs $\mathrm{Sch}(\Gamma /\Gamma _{n},S)$ will
have independence ratio bounded away from $1/2$.

\section{Amenable groups and free groups \label{amensec}}

In this section we discuss weak containment in the realm of amenable groups
and then apply the result for free groups. We also show how to derive a
recent theorem of Conley and Kechris on the maximal measure of independent
subsets for measure preserving actions.

\begin{lemma}
Let $\Gamma $ be an amenable group and let $g$ be a measure preserving
action of $\Gamma $ on a finite set. Then every free measure preserving
ergodic action of $\Gamma $ weakly contains $g$.
\end{lemma}

\noindent \textbf{Proof. }It is known 
(see \cite[13.2]{kechbook} and \cite{forweiss}
), that any two free, measure-preserving actions of an amenable
group are weakly equivalent. Let $b=\{0,1\}^{\Gamma }$ denote the standard
Bernoulli action of $\Gamma $ and let $g_{0}=g\times b$. Then $g_{0}$ is
weakly contained in any measure preserving free action, so the same holds
for its factor $g$. $\square $\bigskip

\begin{lemma}
\label{kislemma}Let $\Gamma $ be a countable group and let $f$ be a measure
preserving action of $\Gamma $. Let $(\Gamma _{n})$ be a chain in $\Gamma $,
let $g_{n}$ be the coset action of $\Gamma $ on $\Gamma /\Gamma _{n}$ and
let $g$ be the boundary action of $\Gamma $ with respect to $(\Gamma _{n})$.
Then $f$ weakly contains $g$ if and only if $f$ weakly contains $g_{n}$ for
all $n$.
\end{lemma}

\noindent \textbf{Proof. }Since $g_{n}$ is a factor of $g$, if $f$ weakly
contains $g$ then it also weakly contains $g_{n}$ for all $n$. In the other
direction, every finite measurable partition of the underlying measure space
of $g$ can be approximated by partitions of the underlying sets of $g_{n}$,
projected to the boundary of the coset tree of $g$ with arbitrarily small
error. Hence, if $f$ weakly contains all the $g_{n}$, then it can simulate
any partition of the underlying set of $g$ as well, and so it weakly
contains $g$. $\square $

\bigskip

\noindent \textbf{Proof of Theorem \ref{frattini}. }Let $F$ be a free group
of rank $d$ and let $p$ be a prime. Let $\mathcal{N}$ be the set of normal
subgroups of $F$ with finite $p$-power index and let $\mathcal{K}$ be the
set of normal subgroups of $F$ where the quotient group is finite and
solvable. For $l>0$ let $\mathcal{K}_{l}\subset \mathcal{K}$ consist of
normal subgroups where the quotient group has derived length at most $l$ and
let $\mathcal{N}_{l}=\mathcal{N\cap K}_{l}$. Let $G$ denote the inverse
limit of $F$ with respect to $\mathcal{N}$; $G$ is called the pro $p$%
-completion of $F$. Let $g$ denote the left action of $F$ on $G$. Since $F$
is residually a $p$-group, $g$ is ergodic and free. Similarly, let $S$
denote the pro (finite solvable) completion of $F$, that is, the inverse
limit of $F$ with respect to $\mathcal{K}$ and let $s$ denote the left
action of $F$ on $S$.

Since every finite $p$-group is solvable, $\mathcal{N}$ is a subset of $%
\mathcal{K}$ and so $g$ is a factor of $s$. In particular, $s$ weakly
contains $g$.

In the other direction, let $h$ be a finite action of $F$ with solvable
image. Let $l$ be the derived length of the image of $h$ and let $\Gamma
=F/F^{(l)}$ be the free solvable group of derived length $l$, where $F^{(l)}$
denotes the $l$-th element of the derived series of $F$. Let $G_{l}$ denote
the inverse limit of $F$ with respect to $\mathcal{N}_{l}$ and let $g_{l}$
denote the left action of $F$ on $G_{l}$. Let $S_{l}$ denote the inverse
limit of $F$ with respect to $\mathcal{K}_{l}$ and let $s_{l}$ denote the
left action of $F$ on $S_{l}$. It is easy to see that $F^{(l)}\leq \mathrm{%
Ker}(s_{l})\leq \mathrm{Ker}(g_{l})$, so in fact $s_{l}$ and $g_{l}$ can 
also be regarded as $\Gamma $-actions. Again, $g_{l}$ is a factor of $s_{l}$.
Now by a result of Gruenberg \cite{gru} $\Gamma $ is residually $p$,
 which implies
that $g_{l}$ (and hence $s_{l}$) are free as $\Gamma $-actions. Since $%
\Gamma $ is amenable, $g_{l}$ weakly contains $s_{l}$ as a $\Gamma $-action.
But then $g_{l}$ weakly contains $s_{l}$ as an $F$-action as well. Now $h$
is a factor for $s_{l}$, so $g_{l}$ weakly contains $h$. But then $g$ weakly
contains $h$ as well, since $g_{l}$ is a factor of $g$. Since $F$ is
finitely generated, it has finitely many subgroups of a given index. Hence $%
\mathcal{K}$ is countable, and so it is generated by a chain in $F$. In
particular, $s$ is a boundary action with respect to a chain in $\mathcal{K}$%
. Using Lemma \ref{kislemma}, $g$ weakly contains $s$.

The theorem is proved. $\square $

\bigskip

\noindent \textbf{Remark.} One can ask whether the whole profinite
completion of a finitely generated free group is weakly equivalent to its
pro $p$ completion. To prove this, it would suffice to show the following:
if $F$ is a finitely generated free group and $N$ is a normal subgroup of
finite index in $F$, then there exists a normal subgroup $K\leq N$ in $F$
such that $F/K$ is amenable and residually $p$.

\bigskip

Now we present how to derive the following recent results of Conley and
Kechris \cite[Theorems 0.5 and 0.6]{conkec}, using the language established
in this paper. For a measure preserving action $a$ of $\Gamma $ on ($X,\mu $%
) and a finite generating set $S$ of $\Gamma $, we call a subset $Y\subseteq
X$ $S$-independent, if for all $y\in Y$ and $s\in S$, $ys\notin Y$. Let $%
i(S,a)$ denote the supremum of $\mu $-measures of $S$-independent Borel
subsets. The same way, we call a $c$-coloring $f:X\rightarrow \{1,\ldots
,c\} $ to be $S$-legal, if for all $x\in X$ and $s\in S$, $f(x)\neq f(xs)$.
Let $\varkappa (S,a)$ denote the minimal $c$ such that $X\ $has an $S$-legal 
$c$-coloring.

\begin{theorem}
Let $\Gamma $ be a countable group and $S\subseteq \Gamma $ a finite
symmetric set of generators with $Cay(\Gamma ,S)$ bipartite. Then the
following are equivalent: \newline
(i) $\Gamma $ is amenable;\newline
(ii) $i(S,a)$ is constant for any free, measure preserving action $a$ of $%
\Gamma $; \newline
(iii) $i(S,a)=1/2$, for any free, measure preserving action $a$ of $\Gamma $%
; \newline
(iv) $\varkappa (S,a)$ is constant for any free, measure preserving action $%
a $ of $\Gamma $; \newline
(v) $\varkappa (S,a)=2$, for every free, measure preserving action $a$ of $%
\Gamma $.
\end{theorem}

\begin{theorem}
Let $\Gamma $ be a countable group and $S\subseteq \Gamma $ a finite
symmetric set of generators with $Cay(\Gamma ,S)$ bipartite. Then the
following are equivalent: \newline
(i)' $\Gamma $ has property (T);\newline
(ii)' $i(S,a)<1/2$, for any free, measure preserving, weakly mixing action $%
a $ of $\Gamma $; \newline
(iii)' $\varkappa (S,a)\geq 3$, for every free, measure preserving, weakly
mixing action $a$ of $\Gamma $.
\end{theorem}

\noindent \textbf{Proof. }Since $Cay(\Gamma ,S)$ is bipartite, $\Gamma $
acts on the two point set such that no element of $S$ fixes a point. Let us
call this action $g$ and its kernel $N$. Let $f$ be the Bernoulli action of $%
\Gamma $ on $\{0,1\}^{\Gamma }$ endowed with the product measure. Then $g$
is not a factor of $f$, since the action of $N$ on $\{0,1\}^{\Gamma }$ is
isomorphic to $\{0,1,2,3\}^{N}$ and hence its ergodic. Let $h$ be the
induced action of the Bernoulli action of $N$ on $\{0,1\}^{N}$ to $\Gamma $.
Then $h$ factors on $g$, so $i(S,h)=1/2$ and $\varkappa (S,h)=2$.

If $\Gamma $ is amenable, then any two free, measure-preserving, ergodic
actions of $\Gamma $ are weakly equivalent. Hence (ii) holds and the
constant has to be $1/2$ by considering $h$. So all of (ii), (iii), (iv) and
(v) holds. If $\Gamma $ is non-amenable, then by \cite{losrind} $f$ is
strongly ergodic, so by Theorem \ref{contrig}, $f$ does not even weakly
contain $g$. In particular, $i(S,f)<1/2$ and $\varkappa (S,f)>2$. Again
considering $h$, we see that all of (ii), (iii), (iv) and (v) fail.

If $\Gamma $ has property (T), then by \cite{schmidt} any free,
measure-preserving, ergodic action $a$ of $\Gamma $ is strongly
ergodic, and by weak mixing, the restriction of $a$ to $N$ stays ergodic, so 
$a$ does not factor on $g$. Hence by Theorem \ref{contrig}, $a$ does not
even weakly contain $g$. In particular, $i(S,a)<1/2$ and $\varkappa
(S,a)\geq 3$. So both (ii)' and (iii)' hold. If $\Gamma $ does not have
property (T), then by \cite{glasnweiss} there exists a free,
measure-preserving, weakly mixing action $a$ of $\Gamma $ that is not
strongly ergodic. By weak mixing, $a$ does not factor on $g$, hence the
restriction of $a$ to $N$ stays ergodic, but not strongly ergodic, and so by
Schmidt's Lemma, it weakly contains $\frac{1}{2}\mathrm{Id}_{N}+\frac{1}{2}%
\mathrm{Id}_{N}$, which is equivalent to saying that $a$ weakly contains $g$%
. So both (ii)' and (iii)' fail.

$\square $

\end{document}